\documentclass[11pt,a4paper]{amsart}
\title[Regularizations of residue currents]{Regularizations of products of residue and principal value currents}

\author{H\aa kan Samuelsson}
\address{Matematik,
  Chalmers tekniska h\"ogskola och G\"oteborgs universitet,\linebreak
  S-412~96  G\"oteborg, Sweden}
\email{hasam@math.chalmers.se}
\date{\today}

\usepackage[english]{babel}
\usepackage{amsmath}
\usepackage{amsthm,amssymb,latexsym}
\usepackage{t1enc}
\usepackage{fancyheadings}
\usepackage{epic}
\usepackage{eepic}
\usepackage{xypic}
\usepackage{mathrsfs}
\usepackage{setspace}

\newtheorem{prop}{Proposition}
\newtheorem{lemma}[prop]{Lemma}
\newtheorem{cor}[prop]{Corollary}
\newtheorem{thm}[prop]{Theorem}

\theoremstyle{definition}

\newtheorem{remark}[prop]{Remark}

\newcommand{\C}{\mathbb{C}}
\newcommand{\D}{\mathscr{D}}
\newcommand{\R}{\mathbb{R}}

\newcommand{\debar}{\bar{\partial}}
\def\newop#1{\expandafter\def\csname #1\endcsname{\mathop{\rm
#1}\nolimits}}

\newop{span}
\usepackage{graphicx}

\begin{document}
\nocite{*}
\bibliographystyle{plain}
\maketitle
\thispagestyle{empty}
 
\begin{abstract}
Let $f_1$ and $f_2$ be two functions on some complex $n$-manifold and let $\varphi$ 
be a 
test form of bidegree $(n,n-2)$. Assume that $(f_1,f_2)$ 
defines a complete intersection. 
The integral of $\varphi/(f_1f_2)$ on 
$\{|f_1|^2=\epsilon_1, |f_2|^2=\epsilon_2\}$ is the 
residue integral
$I_{f_1,f_2}^{\varphi}(\epsilon_1,\epsilon_2)$.
It is in general discontinuous at the origin. Let $\chi_1$ and $\chi_2$ be smooth
functions on $[0,\infty]$ such that $\chi_j(0)=0$ and $\chi_j(\infty)=1$. We prove that 
the regularized residue integral 
defined as the integral of $\debar \chi_1\wedge \debar \chi_2\wedge\varphi/(f_1f_2)$,
where $\chi_j=\chi_j(|f_j|^2/\epsilon_j)$, 
is Hölder continuous on the closed first quarter
and that the value at zero is the Coleff-Herrera residue current acting on $\varphi$.
In fact, we prove that if $\varphi$ is a test form of bidegree $(n,n-1)$ then 
the integral of $\chi_1\debar\chi_2\wedge \varphi/(f_1f_2)$
is Hölder continuous and tends to the $\debar$-potential $[(1/f_1)\wedge \debar (1/f_2)]$
of the Coleff-Herrera current,
acting on $\varphi$.
More generally, let $f_1$ and $f_2$ be sections of some vector bundles
and assume that $f_1\oplus f_2$ defines a complete intersection.
There are associated principal value currents $U^f$ and $U^g$ and residue currents 
$R^f$ and $R^g$. 
The residue currents equal the Coleff-Herrera residue currents locally.
One can
give meaning to formal expressions such as e.g.\ $U^f\wedge R^g$ in such a way that
formal Leibnitz rules hold.
Our results generalize to products of these currents as well.
\end{abstract}

\section{Introduction}
Consider a holomorphic function $f$ defined on some complex $n$-manifold $X$
and let $V_f=f^{-1}(0)$. Schwartz found that there is a distribution, or current, 
$U$ on $X$ 
such that $fU=1$, \cite{Sch2}. 
The existence of the 
principal value current $[1/f]$ defined by 
\[
\D_{n,n}(X) \ni \varphi \mapsto \lim_{\epsilon \rightarrow 0}\int_{|f|^2>\epsilon}\varphi/f
\]
was proved by Herrera and Lieberman in \cite{HL} using 
Hironaka's desingularization theorem, \cite{H} and gives a realization of such a current 
$U$. The $\debar$-image of the principal value
current is the residue current associated to $f$. By Stokes' theorem its action on a
test form of bidegree $(n,n-1)$ is given by the limit as $\epsilon\rightarrow 0$
(along regular values for $|f|^2$) of the residue integral
\begin{equation}\label{hl}
I_f^{\varphi}(\epsilon)=\int_{|f|^2=\epsilon}\varphi/f.
\end{equation}
One main point discovered by Herrera and Lieberman is that if $\varphi$ has bidegree
$(n-1,n)$ then for each $k$, $I_{f^k}^{\varphi}(\epsilon)=\mathcal{O}(\epsilon^{\delta_k})$ 
for some positive
$\delta_k$. Using this, one can then smoothen the integration over 
$|f|^2=\epsilon$ and regularize
the residue current by using smooth functions
$\chi$ defined on $[0,\infty)$ such that $\chi$ is $0$ at zero and tends to $1$
at infinity. In fact, we can make a Leray decomposition and 
write any $(n,n)$-test form $\varphi$ as 
$\phi\wedge\partial f/f^{k}$ for some $k$, where $\phi$ is a test form of bidegree 
$(n-1,n)$ whose restriction to $|f|^2=t$ is unique, for each $t>0$. Then writing the integral 
of $\chi(|f|^2/\epsilon)\varphi/f$ as an integral over the level surfaces $|f|^2=t$ and 
using Herrera's and Lieberman's result one sees that $\chi(|f|^2/\epsilon)/f$ is a 
regularization of the principal value current $[1/f]$. 
It follows that the residue current can be obtained as the weak limit of 
the smooth form $\debar \chi(|f|^2/\epsilon)/f$. This is also a consequence of  
Corollary \ref{spenat} below.
A natural choice for $\chi$ is
$\chi(t)=t/(t+1)$ and we see that we get the well known result that 
the residue current can be obtained as the weak limit 
of $\debar (\bar{f}/(|f|^2+\epsilon))$. 
We also briefly mention the more general currents studied by Barlet, \cite{Barlet}.
If we instead integrate over the fiber $f=s$ in \eqref{hl} and let $\varphi$ have bidegree
$(n-1,n-1)$ then the integral has an asymptotic expansion in $s$ with current coefficients. 
The constant term is Lelong's integration current on $V_f$ 
and the residue current $\debar [1/f]$ can be obtained from the coefficient 
of $s^n$.
  
We turn to the main focus of this paper which is the codimension two case. Let $f$ and 
$g$ be two holomorphic functions on $X$ such that $f$ and $g$ define a complete 
intersection, that is, the common zero set $V_{f\oplus g}$ has codimension two.
Consider the residue integral
\begin{equation}\label{resint}
I_{f,g}^{\varphi}(\epsilon_1,\epsilon_2)=
\int_{\stackrel{\scriptstyle|f|^2=\epsilon_1}{|g|^2=\epsilon_2}}
\frac{\varphi}{fg}.
\end{equation} 
The unrestricted limit of the residue integral 
as $\epsilon_1,\epsilon_2\rightarrow 0$ does not exist in general. The first example of 
this phenomenon was discovered by Passare and Tsikh in \cite{PTmotex}, 
and Björk later found 
that this indeed is the typical case, \cite{jeb1}. See also \cite{Pavlova}. 
Via Hironaka's theorem on resolutions 
of singularities one may assume that the hypersurface $f\cdot g=0$ has normal crossings,
which means that there is a (finite) atlas of charts such that 
$f(\zeta)=\tilde{f}(\zeta)\zeta^{\alpha}$ and 
$g(\zeta)=\tilde{g}(\zeta)\zeta^{\beta}$ where $\alpha$ and
$\beta$ are multiindices (depending on the chart) and 
$\tilde{f}$ and $\tilde{g}$ are invertible holomorphic
functions. It is actually the invertible factors which cause problems. One can always dispose
of one of the factors, but in general not of both. However, if the matrix $A$,
whose two rows are the integer vectors $\alpha$ and $\beta$ respectively, has rank two
there is a change of variables $z=\tau(\zeta)$ such that 
$z^{\alpha}=\tilde{f}(\zeta)\zeta^{\alpha}$ and
$z^{\beta}=\tilde{g}(\zeta)\zeta^{\beta}$, see e.g.\ \cite{krelle}. Hence, when $\alpha$ 
and $\beta$ are not linearly dependent we can make both the invertible factors disappear. 
Problems therefore arise in so called charts of resonance where $\alpha$ and $\beta$ 
are linearly dependent.
Coleff and Herrera realized that if one demands that $\epsilon_1$ and $\epsilon_2$
tend to zero in such a way that $\epsilon_1/\epsilon_2^k\rightarrow 0$ for all 
$k\in \mathbb{Z}_+$, along a so called admissible path, 
then one will get no contributions from the charts of resonance 
because one cannot have $|\tilde{f}(\zeta)\zeta^{\alpha}|<<|\tilde{g}(\zeta)\zeta^{\beta}|$
if $\alpha$ and $\beta$ are linearly dependent.
They proved in \cite{CH} that the limit, along an admissible path, of the residue integral
exists and defines the action of a $(0,2)$-current,
the Coleff-Herrera residue current $[\debar (1/f) \wedge \debar (1/g)]$. 
In \cite{krelle}
Passare smoothened the integration over the set $\{|f|^2=\epsilon_1\}\cap\{|g|^2=\epsilon_2\}$
by introducing functions $\chi$ as described above, and he studied possible weak limits 
of forms 
\begin{equation}\label{pform}
\frac{\debar \chi_1(|f|^2/\epsilon_1)}{f}\wedge\frac{\debar \chi_2(|g|^2/\epsilon_2)}{g}
\end{equation} 
along parabolic paths $(\epsilon_1,\epsilon_2)=(\epsilon^{s_1},\epsilon^{s_2})$ where
$s=(s_1,s_2)$ belongs to the simplex $\Sigma_2(2)=\{(x,y)\in \R_+^2;\, s_1+s_2=2\}$. 
He found that it is enough
to impose finitely many linear conditions $(n_j,s)\neq 0$ 
to assure that \eqref{pform} has a weak limit along the corresponding parabolic path. 
The linear conditions partition $\Sigma_2(2)$ into finitely many open
segments and the weak limit of \eqref{pform} along a parabolic path corresponding to 
an $s$ in such a segment only depends on the segment. We say that $(\epsilon_1,\epsilon_2)$ 
tends to zero inside a 
Passare sector. Moreover, as we assume that $f$
and $g$ define a complete intersection, the limit is even independent of the choice of
segment. In this case it also coincides with the Coleff-Herrera current.
One can obtain a $\debar$-potential to the Coleff-Herrera current e.g.\ by changing 
the integration set in \eqref{resint} to $\{|f|^2>\epsilon_1\}\cap\{|g|^2=\epsilon_2\}$ and 
pass to the limit along an admissible path or by removing the first $\debar$ in 
\eqref{pform} and pass to the limit inside a Passare sector. This $\debar$-potential 
is denoted $[(1/f)\debar (1/g)]$. The main result in this paper implies that if 
$\chi_j\in C^{\infty}([0,\infty])$ satisfy $\chi_j(0)=0$ and $\chi_j(\infty)=1$ then,
in the sense of currents
\begin{equation}\label{Hres}
\lim_{\epsilon_1,\epsilon_2\rightarrow 0}
\frac{\chi_1(|f|^2/\epsilon_1)}{f}\frac{\debar \chi_2(|g|^2/\epsilon_2)}{g}=
\big[\frac{1}{f}\debar\frac{1}{g}\big],
\end{equation}
and the action of the smooth form on the left hand side on a test form depends
Hölder continuously on 
$(\epsilon_1,\epsilon_2)\in [0,\infty)^2$.
For the particular 
case when $\chi_j(t)=t/(t+1)$ our result, apart from the Hölder continuity, 
was announced in \cite{HS}.  
Actually, it is possible to relax the smoothness assumption on one of the $\chi_j$ 
in \eqref{Hres}.
As mentioned above, one can always dispose of one of the invertible factors. Say that we 
always arrange so that $\tilde{f}\equiv 1$. Then, examining the proof, one finds that 
one may take $\chi_1$ to be the characteristic function of $[1,\infty]$.
Hence,
\[
\int_{|f|^2>\epsilon_1}\frac{\debar\chi_2(|g|^2/\epsilon_2)}{fg}\wedge \varphi
\rightarrow 
\big[\frac{1}{f}\debar\frac{1}{g}\big].\varphi
\]
with H\"older continuity.
Note that if we let both $\chi_1$ and $\chi_2$ 
be the characteristic function of $[1,\infty]$ then this result is no longer true 
in view of the examples of Passare-Tsikh and Björk.

Our result also 
generalize to products of pairs of so called Bochner-Martinelli blocks. Consider a tuple
$f=(f_1,\ldots,f_m)$ of holomorphic functions on $X$. The residue integral corresponding
to $f$,
$I_f^{\varphi}(\epsilon_1,\ldots,\epsilon_m)$, is defined analogously to \eqref{resint}.
If we take the mean value of the residue integral over 
$\epsilon=(\epsilon_1,\ldots,\epsilon_m)$ in the simplex 
$\Sigma_m(\delta)=\{s\in \R_+^m;\,\sum s_j=\delta\}$ we obtain
\begin{equation}\label{BM}
c_m\int_{|f|^2=\delta}
\frac{\sum_{j=1}^m (-1)^{j+1}\bar{f}_j\bigwedge_{i\neq j}\debar\bar{f}_i}{|f|^{2m}}
\wedge\varphi,
\end{equation}
where $c_m$ is a constant only depending on $m$. It turns out, see \cite{PTY}, 
that the limit as 
$\delta$ tends to zero of \eqref{BM} exists and defines the action of a $(0,m)$-current,
which in the case $f$ defines a complete intersection, coincides with the Coleff-Herrera
current and also with the currents studied in \cite{BGY} and \cite{PTcanada}.
Based on the work in \cite{PTY} Andersson introduces more general currents of the 
Cauchy-Fantappiè-Leray type in \cite{A1}. We will briefly discuss Andersson's construction
in Section \ref{prel&not}. In short, he defines a singular form $u^f=\sum u^f_{k,k-1}$, 
where the terms $u^f_{k,k-1}$ are similar to the form in \eqref{BM}, 
and he shows that it is extendible to $X$ as a current,
$U^f$, either as principal values or by analytic continuation. The residue current,
$R^f$, is derived from the current $U^f$ and equals the Coleff-Herrera current 
locally if $f$ defines a complete intersection. If $g$ is also a tuple of functions 
there is a natural way of defining the product of the Cauchy-Fantappiè-Leray type 
currents corresponding to $f$ and $g$ so that formal Leibnitz rules hold, see \cite{W}.
If $f\oplus g$ defines a complete intersection and 
$\chi_1,\chi_2\in C^{\infty}([0,\infty])$ vanish to high enough orders at zero and 
equals $1$ at infinity then we prove that the smooth forms
\[
\chi_1(|f|^2/\epsilon_1)u^f\wedge \debar \chi_2(|g|^2/\epsilon_2)\wedge u^g\,\,\,\mbox{and}
\]
\[ 
\debar\chi_1(|f|^2/\epsilon_1)\wedge u^f\wedge \debar \chi_2(|g|^2/\epsilon_2)\wedge u^g
\]
are Hölder continuous as currents for $(\epsilon_1,\epsilon_2)\in [0,\infty)^2$
and tend to $U^f\wedge R^g$ and $R^{f}\wedge R^g$ respectively as 
$\epsilon_1,\epsilon_2\rightarrow 0$; Theorem \ref{satsen} and Corollary \ref{sven}.
If $g$ is a function such that $f\oplus g$ defines a complete intersection, 
our techniques can also be used to prove that 
$\debar\chi_1(|f|^2/\epsilon_1)\wedge u^f\chi_2(|g|^2/\epsilon_2)\rightarrow R^f$
when $\chi_2$ equals the characteristic function of $[1,\infty]$.
We use this to conclude that $R^f$ has the 
standard extension property in the complete intersection case, Corollary \ref{SEP}.
For more historical accounts we refer to the survey article \cite{jeb2} by Björk.

The disposition of the paper is as follows: In Section \ref{schysstsekt} 
we outline a proof of 
\eqref{Hres} since the proofs of the more general statements about Bochner-Martinelli
or Cauchy-Fantappiè-Leray blocks are only more difficult to prove in the technical sense
and to make it clear that it is not necessary to work through the constructions
of Bochner-Martinelli or Cauchy-Fantappiè-Leray type currents in order to prove 
\eqref{Hres}. In Section \ref{prel&not} we recall Andersson's construction and explain 
some useful notation. Section \ref{CFLsektion} 
contains some fairly well known regularization 
results about Cauchy-Fantappiè-Leray type currents. As Andersson's formalism makes the 
arguments a little smoother we also supply the proofs. Section \ref{monom}
contains the technical core of this paper. We study regularizations of products of monomial
currents which we then use in Section \ref{prodsektion} to prove our main results;
Theorem \ref{satsen} and its corollaries \ref{sven}, \ref{leif}
and \ref{CR} and Theorem \ref{kent}. In Section \ref{motexsekt}
we see by explicit computations that Corollary \ref{CR} holds for the example by
Passare and Tsikh. This section is essentially self-contained.  

\section{Sketch of proof in the case of two functions}\label{schysstsekt}
Let $f$ and $g$ be two holomorphic functions on $X$ defining a complete intersection. 
We sketch how one can handle the difficulties arising in charts of resonance when
proving \eqref{Hres}. We study the integral
\begin{equation}\label{janne}
\int \frac{\chi_1(|f|^2/\epsilon_1)}{f}\frac{\debar \chi_2(|g|^2/\epsilon_2)}{g}
\wedge \varphi
\end{equation}
where $\varphi$ is a test form of bidegree $(n,n-1)$. By Hironaka's theorem we may assume 
that $f=\zeta^{\alpha}\tilde{f}$ and 
$g=\zeta^{\beta}\tilde{g}$ are
monomials times non-vanishing functions. One of the non-zero factors can be incorporated 
in a variable and so we assume that $\tilde{f}\equiv 1$.
We assume also that we are in a chart of 
resonance, i.e.\ that $\alpha$ and $\beta$ are linearly dependent. After resolving 
singularities $f$ and $g$ no longer define a complete intersection in general,
but on the other 
hand a degree argument shows that  
$d\bar{\zeta}_j/\bar{\zeta}_j\wedge \varphi$ becomes a test form for any $\zeta_j$ 
dividing both $f$ and $g$. See the proof of Theorem \ref{satsen} for more details. 
Since $\alpha$ and $\beta$ are linearly dependent, 
$d\bar{\zeta}_j/\bar{\zeta}_j\wedge \varphi$ is a test form for all $j$ such that 
$\alpha_j\neq 0$, or equivalently, $\beta_j\neq 0$.
Now, \eqref{janne} equals
\[
\sum_j \beta_j\int 
\frac{\chi_1(|\zeta^{\alpha}|^2/\epsilon_1)}{\zeta^{\alpha}}
\frac{\chi'_2(\Psi|\zeta^{\beta}|^2/\epsilon_2)}{\zeta^{\beta}}
\frac{|\zeta^{\beta}|^2}{\epsilon_2}\wedge
\frac{d\bar{\zeta}_j}{\bar{\zeta_j}}
\wedge \varphi/\tilde{f}
\]
where $\Psi=|\tilde{g}|^2$ is a strictly positive smooth function. It now follows from
Corollary \ref{korret} that each term in this sum tends to zero as $\epsilon_1$ and 
$\epsilon_2$ tend to zero.
Hence the charts of resonance do not give any contributions.

\section{Preliminaries and notation}\label{prel&not}
Assume that $f$ is a section of the dual bundle $E^*$ of a holomorphic 
$m$-bundle $E\rightarrow X$ over a complex $n$-manifold $X$.
We will only deal with local problems and it is therefore no loss of generality in
assuming that $E\rightarrow X$ is trivial. However, the formalism will run 
smoother with an invariant notation. As mentioned above, we will recall Andersson's 
construction in \cite{A1} and produce currents $U^f$ and $R^f$ and we emphasize that
in the case $E\rightarrow X$ is the trivial line bundle then $U^f$ and $R^f$ are the 
currents $[1/f]$ and $\debar[1/f]$ times some basis elements.
On the exterior algebra $\Lambda E$ of $E$, the section $f$
induces mappings $\delta_f\colon \Lambda^{k+1}E\rightarrow \Lambda^kE$ 
of interior multiplication and $\delta_f^2=0$. We introduce the spaces 
$\mathcal{E}_{0,q}(X,\Lambda^kE)$ of the smooth sections of the exterior 
algebra of $E\oplus T^*_{0,1}X$ which are $(0,q)$-forms with values in
$\Lambda^kE$. We also introduce the corresponding spaces of currents,
$\D'_{0,q}(X,\Lambda^kE)$. The mappings $\delta_f$ extend to mappings
$\delta_f\colon \D'_{0,q}(X,\Lambda^{k+1}E)\rightarrow\D'_{0,q}(X,\Lambda^kE)$
with $\delta_f^2=0$ and these mappings anti-commute with the $\debar$-operator.
Hence, $\D'_{0,q}(X,\Lambda^kE)$ is a double complex and the associated 
total complex is
\[
\cdots \stackrel{\nabla_f}{\rightarrow} \mathcal{L}^{r-1}(X,E)
\stackrel{\nabla_f}{\rightarrow} 
 \mathcal{L}^{r}(X,E)\stackrel{\nabla_f}{\rightarrow}\cdots
\]
where $\mathcal{L}^{r}(X,E)=\bigoplus_{q-k=r}\D'_{0,q}(X,\Lambda^kE)$ and
$\nabla_f=\delta_f-\debar$. We will refer to the total complex as the Andersson complex.
The exterior product, $\wedge$, induces mappings
\[
\wedge\colon \mathcal{L}^{r}(X,E)\times \mathcal{L}^{s}(X,E)\rightarrow
\mathcal{L}^{r+s}(X,E)
\]
when possible, and $\nabla_f$ is an antiderivation, i.e.\ $\nabla_f(\tau\wedge \sigma)=
\nabla_f\tau\wedge \sigma+(-1)^r\tau\wedge \nabla_f\sigma$ if 
$\tau\in \mathcal{L}^{r}(X,E)$ and $\sigma\in \mathcal{L}^{s}(X,E)$.
If 
$\tau\in\mathcal{L}^{r}(X,E)$ we write $\tau_{k,k+r}$ for the component of 
$\tau$ belonging to $\D'_{0,k+r}(X,\Lambda^kE)$. Note that functions define 
elements of $\mathcal{L}^{0}(X,E)$ of degree $(0,0)$ and sections of $E$ 
define elements of $\mathcal{L}^{-1}(X,E)$ of degree $(1,0)$. One can show,
see \cite{A1}, that 
if $X$ is Stein and the zero:th cohomology group of the Andersson complex
vanishes then for any holomorphic function $h$ there is a holomorphic 
section $\psi$ of $E$ such that $\delta_f\psi=h$. This means that if 
$f=(f_1,\ldots,f_m)$ in some local holomorphic frame for $E^*$ then 
the division problem $\sum f_j\psi_j=h$ has a holomorphic solution. 
This cannot hold for all $h$ if $f$ has zeros and the Andersson complex can 
therefore not be exact in this case. Still, we try to look for an element 
$u^f\in \mathcal{L}^{-1}(X,E)$ such that $\nabla_fu^f=1$. To this end we assume 
that $E$ is equipped with some Hermitian metric $|\cdot |$ and we let $s_f$
be the section of $E$ with pointwise minimal norm such that $\delta_fs_f=|f|^2$.
Outside $V_f=f^{-1}(0)$ we may put
\[
u^f=\frac{s_f}{\nabla_fs_f}
=\frac{s_f}{\delta_fs_f-\debar s_f}=
\sum_k\frac{s_f\wedge (\debar s_f)^{k-1}}{|f|^{2k}}.
\]
Observe that $\nabla_fs_f$ has even degree so the expression 
$s_f/\nabla_fs_f$ has meaning outside $V_f$ and it follows immediately
that $\nabla_fu=1$ there. The following theorem is proved in \cite{A1}.
\begin{thm}
Assume that $f$ is locally nontrivial. The forms $|f|^{2\lambda}u^f$ and
$\debar |f|^{2\lambda}\wedge u^f$ are locally bounded if $\mathfrak{Re}\,\lambda$
is sufficiently large and they have analytic continuations as currents to 
$\mathfrak{Re}\,\lambda > -\epsilon$. Let
$U^f$ and $R^f$ denote the values at $\lambda=0$. Then $U^f$ is a current 
extension of $u^f$, $R^f$ has support on $V_f$ and 
\[
\nabla_f U^f=1-R^f.
\]
Moreover, $R^f=R^f_{p,p}+\cdots+R^f_{q,q}$
where $p=\mbox{Codim}(V_f)$ and $q=\mbox{min}(m,n)$.
\end{thm}
Note that if $V_f=\emptyset$ then $\nabla_f U^f=1$ on all of $X$, 
which implies that
taking the exterior product with $U^f$ is a homotopy operator for the Andersson
complex. The current $R^f$ is the Bochner-Martinelli, or more generally, 
the Cauchy-Fantappiè-Leray current associated to $f$, and if $f=(f_1,\ldots,f_m)$ in some 
local holomorphic frame, $e_1,\ldots,e_m$, of $E$ then 
\begin{equation}\label{BM=CH}
R^f=\big[\debar\frac{1}{f_1}\wedge\cdots\wedge\debar\frac{1}{f_m}\big]\wedge
e_1\wedge\cdots\wedge e_m
\end{equation}
if $f$ defines a complete intersection, see \cite{A1}.

Now if $f_j$, $j=1,2$,
are sections of the dual bundles $E_j^*$ of holomorphic Hermitian 
$m_j$-bundles $E_j\rightarrow X$ we can apply the above construction to 
the section $f=f_1\oplus f_2$ of the bundle
$E_1^*\oplus E_2^*$ and obtain the currents
$U^f$ and $R^{f}$. We could also try to combine the individual currents
$U^{f_j}$ and $R^{f_j}$. It is shown in \cite{W} that the forms
\[
|f_1|^{2\lambda}u^{f_1}\wedge|f_2|^{2\lambda}u^{f_2},\,\,\,
|f_1|^{2\lambda}u^{f_1}\wedge\debar |f_2|^{2\lambda}\wedge u^{f_2}\,\,\, \mbox{and} \,\,\,
\debar |f_1|^{2\lambda}\wedge u^{f_1}\wedge\debar |f_2|^{2\lambda}\wedge u^{f_2},
\]
which are locally bounded if $\mathfrak{Re}\, \lambda $ is large enough, have
current extensions to $\mathfrak{Re}\, \lambda >-\epsilon$. The values at 
$\lambda=0$ are denoted 
$U^{f_1}\wedge U^{f_2}$, $U^{f_1}\wedge R^{f_2}$, and $R^{f_1}\wedge R^{f_2}$, 
respectively, 
and formal computation rules such as e.g.\ $\nabla_f (U^{f_1}\wedge R^{f_2})=
(1-R^{f_1})\wedge R^{f_2}=R^{f_2}-R^{f_1}\wedge R^{f_2}$ hold.
It is also shown in \cite{W} that if $f$ defines a complete intersection
then $R^f=R^{f_1}\wedge R^{f_2}$. 

We will use the names $f$ and $g$, rather then $f_1$ and $f_2$, for the 
sections of the two bundles and the symbol $\nabla$, without subscript, 
always denotes $\nabla_{f\oplus g}$.  
We will use multiindices extensively in the sequel. Multiindices will be denoted
$\alpha$ and $\beta$ or $I$ and $J$ and sometimes also $r$ and $\rho$. 
The number of variables will always be 
$n$ but it will be convenient to define multiindices by expressions like
$\alpha=(\alpha_j)_{j\in K}$ for $K\subseteq \{1,\ldots,n\}$.
By this we mean that $\alpha=(a_1,\ldots,a_n)$ where $a_j=0$ if 
$j\notin K$ and $a_j=\alpha_j$ if $j\in K$. Hence, if $z=(z_1,\ldots,z_n)$
then $z^{\alpha}=\prod_{j\in K}z_j^{\alpha_j}$ and similarly for 
$\partial^{\alpha} /\partial z^{\alpha}$. Multiindices are added and multiplied
by numbers as elements in $\mathbb{Z}^n$ and 
$\alpha \pm 1=(\alpha_1\pm 1,\ldots, \alpha_n\pm 1)$. Also, $|\alpha|$ 
denotes the length of $\alpha$ as a vector in Euclidean space and $\# \alpha$
is the cardinality of the support of $\alpha$. 

Integration over domains in
$\C^n$ will always be with respect to the volume form 
$(i/2)^ndz_1\wedge d\bar{z}_1\wedge \ldots \wedge dz_n\wedge d\bar{z}_n:=
(i/2)^ndz\wedge d\bar{z}$ if nothing else
is said. If $\Delta$ is a Reinhardt domain in $\C^n$ and 
$\varphi$ is a function which only depends on the moduli of 
the variables and such that $z^{\alpha}\varphi(z)$ is integrable on $\Delta$
then
\[
\int_{\Delta}z^{\alpha}\varphi(z)=0
\]
if $\alpha$ is a non-zero multiindex. 
This simple fact will play a fundamental role to us 
in what follows and we will refer to it as {\em anti-symmetry}.

Unless otherwise stated, the symbol $\chi$ with various subscripts will 
always denote a smooth function
on $[0,\infty]$ which is zero to some order at $0$ and such that 
$\chi(\infty)=1$. By smooth at infinity we mean that $t\mapsto \chi(1/t)$ is smooth
at zero.

\section{Regularizations of Cauchy-Fantappi\`e-Leray type currents}\label{CFLsektion}
Consider a function $\chi$ as above and 
let $\tilde{\chi}(s)=\chi(1/s)$. Then $\tilde{\chi}$ is differentiable at
$s=0$ and $\tilde{\chi}'(s)=-\chi'(1/s)/s^2$. Letting $t=1/s$ we see that
$\chi'(t)=\mathcal{O}(1/t^2)$ as $t \rightarrow \infty$. This simple 
observation will be frequently used in the sequel. It follows that for any 
continuous function $\varphi$ with compact support in $[0,\infty)$ we have 
$|\varphi(\epsilon t) \chi'(t)|\leq C(t+1)^{-2}$ for a constant independent 
of $\epsilon$. Hence by the dominated convergence theorem we see that
\[
\int_0^{\infty}\frac{d}{dt}\chi(t/\epsilon) \varphi(t) dt=
\int_0^{\infty}\frac{d}{d\tau}\chi(\tau) \varphi(\epsilon\tau) d\tau 
\rightarrow
\varphi(0)\int_0^{\infty}\frac{d}{d\tau}\chi(\tau)d\tau =\varphi(0),
\]
and we have proved
\begin{lemma}\label{apid}
Let $\chi \in C^1([0,\infty])$ satisfy $\chi(0)=0$ and $\chi(\infty)=1$. 
Then $(d/dt) \,\chi (t/\epsilon)\rightarrow 
\delta_0$ as measures on $[0,\infty)$.
\end{lemma}

\begin{prop}\label{fundpropCFL}
Assume $\chi \in C^{\infty}([0,\infty])$ vanishes to order $\ell$ at $0$ and 
satisfies $\chi(\infty)=1$. Then
\[
\lim_{\epsilon \rightarrow 0^+}\int
\chi(|f|^2/\epsilon)u^f_{\ell,\ell-1}\wedge \varphi =
U^f_{\ell,\ell-1}.\varphi
\]
for any test form $\varphi$. 
\end{prop}

\begin{proof}
On the set $\Omega=\{(z,t)\in \C^n\times (0,\infty); |f(z)|^2>t\}$ 
we have, for all fixed $\epsilon >0$, that 
\begin{eqnarray*}
\big|u^f_{\ell,\ell-1}
\frac{d}{dt}\chi(t/\epsilon)\wedge \varphi \big|
& \leq &
C\frac{1}{|f|^{2\ell-1}}\big|\frac{d}{dt}\chi(t/\epsilon)\big|\leq \\
& &
C\frac{t^{1/2}}{t^\ell}\big|\frac{d}{dt}\chi(t/\epsilon)\big|\leq
C\frac{1}{t^{1/2}}\nonumber
\end{eqnarray*}
since $\frac{d}{dt}\chi(t/\epsilon)=\mathcal{O}(t^{\ell-1})$. Hence we have an 
integrable singularity on $\Omega$ and by Fubini's theorem we get
\begin{eqnarray}\label{eq01}
\int_0^{\infty}\frac{d}{dt}\chi(t/\epsilon)
\int_{|f|^2>t}u^f_{\ell,\ell-1}\wedge \varphi \,dt
&=&
\int u^f_{\ell,\ell-1}\wedge \varphi
\int_0^{|f|^2}\frac{d}{dt}\chi(t/\epsilon)dt=\nonumber \\
& &
\int u^f_{\ell,\ell-1}
\chi(|f|^2/\epsilon)\wedge \varphi.
\end{eqnarray}
But $J(t)=\int_{|f|^2>t}u^f_{\ell,\ell-1}\wedge \varphi$ is a continuous function
with compact support in $[0,\infty)$ with $J(0)=U^f_{\ell,\ell-1}.\varphi$, see
\cite{PTY} or \cite{A1}. 
Hence by Lemma \ref{apid} the left hand side of \eqref{eq01}
tends to $U^f_{\ell,\ell-1}.\varphi$ and the proof is complete. 
\end{proof}

If we take $\chi(t)$ equal to appropriate powers of $t/(t+1)$ 
we obtain the following natural ways to 
regularize the currents $U^f$ and $R^f$.
\begin{cor}\label{sture} 
For any test form $\varphi$ we have
\begin{equation}\label{eq02}
\lim_{\epsilon\rightarrow 0^+}\int \sum_{\ell\geq 1}
\frac{s_f\wedge (\debar s_f)^{\ell-1}}{(|f|^2+\epsilon)^\ell}\wedge \varphi =
U^f.\varphi 
\end{equation}
and
\begin{equation}\label{eq03}
\lim_{\epsilon\rightarrow 0^+}\int \sum_{\ell\geq 1} \epsilon
\frac{(\debar s_f)^{\ell}}{(|f|^2+\epsilon)^{\ell+1}}\wedge \varphi =
R^f.\varphi.
\end{equation}
\end{cor}

\begin{proof}
Letting $\chi_\ell(t)=t^{\ell}/(t+1)^{\ell}$ 
we see that 
\[
u^f_{\ell,\ell-1}\chi_\ell(|f|^2/\epsilon)=
\frac{s_f\wedge (\debar s_f)^{\ell-1}}{(|f|^2+\epsilon)^\ell}
\]
and so \eqref{eq02} follows from Proposition \ref{fundpropCFL}.
To show that \eqref{eq03} holds we first note that
\[
\sum_{\ell\geq 1}
\frac{s_f\wedge (\debar s_f)^{\ell-1}}{(|f|^2+\epsilon)^\ell}=
\frac{s_f}{\nabla_f s_f+\epsilon}.
\]
Hence
\[
\nabla_f \sum_{\ell\geq 1}
\frac{s_f\wedge (\debar s_f)^{\ell-1}}{(|f|^2+\epsilon)^\ell}=
\nabla_f \frac{s_f}{\nabla_f s_f+\epsilon}=
\frac{\nabla_f s_f}{\nabla_f s_f+\epsilon}=
1-\sum_{\ell\geq 0} \epsilon
\frac{(\debar s_f)^{\ell}}{(|f|^2+\epsilon)^{\ell+1}}.
\]
Since differentiation is a continuous operation on distributions
it follows from \eqref{eq02} that
\[
\lim_{\epsilon \rightarrow 0^+}
 1-\sum_{\ell\geq 0} \epsilon
\frac{(\debar s_f)^{\ell}}{(|f|^2+\epsilon)^{\ell+1}}=
\nabla_f \lim_{\epsilon \rightarrow 0^+}\sum_{\ell\geq 1}
\frac{s_f\wedge (\debar s_f)^{\ell-1}}{(|f|^2+\epsilon)^\ell}=
\nabla_f U^f = 1-R^f
\]
in the sense of currents. The term with $\ell=0$ in the sum on the left is easily
seen to tend to zero in the sense of currents and hence \eqref{eq03}
follows.
\end{proof}

Note that it is the difference 
\begin{equation}\label{rut}
\debar(\chi_{\ell}u^f_{\ell,\ell-1})-
\delta_f(\chi_{\ell+1}u^f_{\ell+1,\ell})=
\debar \chi_{\ell}\wedge u^f_{\ell,\ell-1}+
(\chi_{\ell}-\chi_{\ell+1})\delta_fu^f_{\ell+1,\ell}
\end{equation}
which converges to the term of $R^f$ of bidegree $(\ell,\ell)$. It is only for the term
of top degree, the last term in \eqref{rut} is not present. This explains why 
the regularization result in \cite{PTY}, Theorem $2.1$, coincides with our result for 
the top degree term but not for the terms of lower degree.

We can also take one $\chi$ which vanishes to high enough order at zero 
to regularize all terms of $U^f$ and $R^f$.

\begin{cor}\label{spenat}
Assume that $\chi \in C^{\infty}([0,\infty])$, vanishes to order 
$\mbox{min}(m,n)+1$ at zero and satisfies $\chi(\infty)=1$. Then for any 
test form $\varphi$ we have
\begin{equation}\label{eq04}
\lim_{\epsilon\rightarrow 0^+}\int \chi(|f|^2/\epsilon)u^f\wedge \varphi=
U^f.\varphi
\end{equation}
\begin{equation}\label{eq05}
\lim_{\epsilon\rightarrow 0^+}\int
\debar\chi(|f|^2/\epsilon)\wedge u^f\wedge \varphi=
R^f.\varphi.
\end{equation}
\end{cor}

\begin{proof}
The first statement follows immediately from Proposition \ref{fundpropCFL}.
For the second one we note that
\[
\nabla \chi u^f=\nabla \chi \wedge u^f + \chi \nabla u^f=
-\debar \chi \wedge u^f +\chi \nabla u^f,
\]
and since $\chi$ vanishes to high enough order at zero all terms are 
smooth. Outside $\{f=0\}$ we have $\nabla u^f=1$ and hence 
$\chi \nabla u^f=\chi$ everywhere. Moreover, $\chi(|f|^2/\epsilon)$ tends
to $1$ in the sense of currents and hence
\[
 \debar \chi \wedge u^f= \chi \nabla u^f-\nabla \chi u^f \rightarrow
1-(1-R^f)=R^f
\]
in the sense of currents.
\end{proof}


\section{Regularizations of products of monomial currents}\label{monom}
This section contains the technical result about the normal crossing case needed 
to prove our main theorems in the next section. Of particular importance is
Proposition \ref{tekniskprop}.
First we need a generalization of 
Taylor's formula. Lemma \ref{taylorlemma} enables us to approximate a smooth
function defined on $\C^n$ in a neighborhood of the union of the coordinate hyperplanes
instead of in a neighborhood of their intersection as in the usual Taylor's formula.
The approximating functions are in our case not polynomials in general 
but have enough similarities for our purposes.
For tensor products of one-variable functions this corresponds to 
multiplying the individual Taylor expansions. 
Lemma \ref{taylorlemma} appears as Lemma $2.3$ in \cite{HS} but the formulation there is 
unfortunately not completely correct.
We also remark that Lemma \ref{taylorlemma}
is very similar to Lemma $2.4$ in \cite{CH} and that very general Taylor expansions are 
considered in Chapter $1$ in \cite{LH1}.
Define the linear operator $M_j^{r_j}$ on $C^{\infty}(\C^n)$ to be the 
operator that maps $\varphi$ to the Taylor polynomial of degree $r_j$ of the
function $\zeta_j \mapsto \varphi(\zeta)$ (centered at $\zeta_j=0$). We note 
that $M_j^{r_j}$ and 
$M_i^{r_i}$ commute. To see this we only need to observe 
that 
\[
\frac{\partial}{\partial \tilde{\zeta_i}}\big(\frac{\partial \varphi}{\partial 
\tilde{\zeta_j}}\big|_{\zeta_j=0}\big)\big|_{\zeta_i=0}=
\frac{\partial^2 \varphi}{\partial \tilde{\zeta_i} \partial \tilde{\zeta_j}}\big|_{
\zeta_i=\zeta_j=0}=
\frac{\partial}{\partial \tilde{\zeta_j}}\big(\frac{\partial \varphi}{\partial 
\tilde{\zeta_i}}\big|_{\zeta_i=0}\big)\big|_{\zeta_j=0}
\]
where $\partial/\partial \tilde{\zeta_j}$ means that we do not specify
whether we differentiate with respect to $\zeta_j$ or $\bar{\zeta_j}$.

\begin{lemma}\label{taylorlemma}
Let $K \subseteq \{1,\ldots,n\}$ have cardinality $\kappa$ and let 
$r=(r_j)_{j\in K}$.
Define the linear operator $M^r_K$ on $C^{\infty}(\C^n)$ by
\[
M^r_K=\sum_{j \in K}M^{r_j}_j-\sum_{\stackrel{i,j\in K}{i<j}}M_i^{r_i}M_j^{r_j}
+\cdots + (-1)^{\kappa+1}M_{j_1}^{r_{j_1}}\cdots M_{j_{\kappa}}^{r_{j_{\kappa}}}.
\]
Then for any $\varphi \in C^{\infty}(\C^n)$ we have
\begin{equation}\label{tayloreq}
\varphi(\zeta) = M^r_K \varphi(\zeta) + \int_{[0,1]^{\kappa}}\frac{(1-t)^{r}}{r!}
\frac{\partial^{r+1}}{\partial t^{r+1}}
\varphi(t\zeta) \,dt
\end{equation}
where $t\zeta$ should be interpreted as $(\xi_1,\ldots,\xi_n)$, $\xi_j=t_j
\zeta_j$ if $j \in K$ and $\xi_j=\zeta_j$ if $j \notin K$. In particular
$\varphi-M^r_K\varphi=\mathcal{O}(|\zeta^{r+1}|)$.
Moreover, $M^r_K \varphi$ can be written as a finite sum of terms, 
$\varphi_{IJ}(\zeta)\zeta^I\bar{\zeta}^J$, with the following properties:
\begin{itemize}
\item[(a)] $\varphi_{IJ}(\zeta)$ is independent of some variable and in particular 
of variable $\zeta_j$ if $I_j+J_j>0$, 
\item[(b)] $I_j+J_j\leq r_j$ for $j\in K$,
\item[(c)] if $L$ is the set of indices $j\in K$ such that 
$\zeta_j \mapsto \varphi_{IJ}(\zeta)$ is non-constant then
$\varphi_{IJ}(\zeta)=\mathcal{O}(\prod_{j \in L}|\zeta_j|^{r_j+1})$.
\end{itemize}
\end{lemma}

\begin{proof}
It is enough to prove the lemma when $K=\{1,\ldots, n\}$.
In case $n=1$,  \eqref{tayloreq} is Taylor's formula. For $n\geq 2$,
we write the integral in \eqref{tayloreq} as an iterated integral. Formula \eqref{tayloreq}
then follows by induction.
One can also show \eqref{tayloreq} by repeated 
integrations by parts. 
The difference $\varphi-M^r_K\varphi$ 
is seen to be of the desired size after performing the differentiations
of $\varphi(t\zeta)$ with respect to $t$ inside the integral. 
To see that $M_K^r\varphi$ can be 
written as a sum of terms $\varphi_{IJ}(\zeta)\zeta^I\bar{\zeta}^J$
with the properties (a), (b), and (c), we let
$r_{\tilde{K}}$, for any $\tilde{K}\subseteq K$, denote the multiindex
$(r_{j_1},\ldots,r_{j_{|\tilde{K}|}})$, $r_{i_j}\in \tilde{K}$. A straight
forward computation now shows that
\begin{eqnarray*}
M_K^r\varphi &=& \sum_{j\in K}M_j^{r_j}(\varphi-M_{K\setminus \{j\}}^{r_{
K\setminus \{j\}}}\varphi) \\
&+&
\sum_{\stackrel{i,j\in K}{i<j}}M_i^{r_i}M_j^{r_j}(\varphi-
M_{K\setminus \{i,j\}}^{r_{K\setminus \{i,j\}}}\varphi) \\
&\vdots& \\
&+&
M_{j_1}^{r_{j_1}}\cdots 
M_{j_{\kappa}}^{r_{j_{\kappa}}}\varphi.
\end{eqnarray*}
From the first part of the proof (and the definition of $M_j^{r_j}$) it 
follows that every term on the right hand side is a finite sum of terms with
the stated properties. 
\end{proof}

\begin{lemma}\label{randtermslemma}
Let $\alpha$ be a multiindex and let $M=M^r_K$ be the operator defined in Lemma 
\ref{taylorlemma} with 
$K$ the set of indices $j$ such that $\alpha_j \geq 2$ and $r_j=\alpha_j-2$, $j\in K$. 
Then for any $\varphi \in \D(\C^n)$ we have
\[
\int_{\Delta}\frac{1}{\zeta^{\alpha}}(\varphi - M\varphi) =
\Big[\frac{1}{\zeta^{\alpha}}\Big]. \, \varphi \, (i/2)^n d\zeta \wedge d\bar{\zeta} 
\]
if $\Delta$ is a polydisc containing the support of $\varphi$.
\end{lemma}

\begin{proof}
Note that by Lemma \ref{taylorlemma} we have $\varphi - M\varphi=\mathcal{O}(
|\zeta^{\alpha-1}|)$ and so $(\varphi - M\varphi)/\zeta^{\alpha}$ 
is integrable on $\Delta$. Hence if we let $\Delta_{\delta}=\Delta \cap_{
j}\{|\zeta_j|>\delta\}$ we get
\begin{eqnarray*}
\int_{\Delta}\frac{1}{\zeta^{\alpha}}(\varphi - M\varphi) &=&
\lim_{\delta \rightarrow 0}
\int_{\Delta_{\delta}}\frac{1}{\zeta^{\alpha}}(\varphi - M\varphi) \\
&=&
\lim_{\delta \rightarrow 0}
\int_{\Delta_{\delta}}\frac{1}{\zeta^{\alpha}}\varphi-
\lim_{\delta \rightarrow 0}\
\int_{\Delta_{\delta}}\frac{1}{\zeta^{\alpha}}M\varphi.
\end{eqnarray*}
The first limit on the right hand side is the tensor product of the principal value 
currents $[1/\zeta_j^{\alpha_j}]$ (acting on $\varphi \,(i/2)^nd\zeta\wedge d\bar{\zeta}$)
and hence it equals $[1/\zeta^{\alpha}].\varphi \,(i/2)^nd\zeta\wedge d\bar{\zeta}$.
It follows by anti-symmetry that actually
\[
\int_{\Delta_{\delta}}\frac{1}{\zeta^{\alpha}}M\varphi=0
\]
for all $\delta >0$. In fact, $M\varphi$
is a sum of terms $\varphi_{IJ}(\zeta)\zeta^I\bar{\zeta}^J$ where $I_j+J_j \leq
\alpha_j-2$ for all $j$ and the coefficient $\varphi_{IJ}(\zeta)$ is at least 
independent of
some variable. 
\end{proof}

\begin{lemma}\label{trollerilemma}
Let $\chi_1,\chi_2\in C^{\infty}([0,\infty])$ and let $\Phi$ and $\Psi$ be 
smooth strictly positive functions on $\C^n$. Let also $M_K^r$ be the 
operator defined in Lemma \ref{taylorlemma} with $K$ and $r$ arbitrary. Then
\[
\chi_1(t_1\Phi)\chi_2(t_2\Psi) =
M_K^r(\chi_1(t_1\Phi)\chi_2(t_2\Psi))
+
|\zeta^{r+1}|B(t_1, t_2, \zeta),
\]
where $B$ is bounded on $(0,\infty)^2\times D$ if $D\Subset \C^n$.
\end{lemma}

\begin{proof}
If $D\Subset \C^n$ both $\Phi$ and $\Psi$ have strictly positive infima and 
finite suprema on $D$ and so there is a neighborhood $U$ of $[0,\infty]^2$
in $\widehat{\R}\times\widehat{\R}$ such that the function
$(t_1,t_2,\zeta)\mapsto \chi_1(t_1\Phi)\chi_2(t_2\Psi)$ is smooth 
on $U\times D$. From Lemma \ref{taylorlemma} it follows that
\[
\chi_1(t_1\Phi)\chi_2(t_2\Psi)=
M_K^r(\chi_1(t_1\Phi)\chi_2(t_2\Psi))+
\sum\limits_{\stackrel{\scriptstyle I,J\subseteq K}{I_j+J_j=r_j+1}}
G_{IJ}(t_1,t_2,\zeta)\zeta^I\bar{\zeta}^J
\]
for some functions $G_{IJ}$ which are smooth on $U\times D$, and the lemma 
readily follows.
\end{proof}

To prove Proposition \ref{tekniskprop} we will need the 
estimates of the following two elementary lemmas.

\begin{lemma}\label{uppsk_lemma2}
Let $\Delta$ be the unit polydisc in $\C^n$ and put 
$\Delta_{\epsilon}^{\alpha}=\{\zeta\in\Delta; |\zeta^{\alpha}|^2\geq \epsilon\}$ and
$\Delta_{\epsilon_1,\epsilon_2}^{\alpha, \beta}=\{\zeta\in\Delta; |\zeta^{\alpha}|^2\geq \epsilon_1,|\zeta^{\beta}|^2\geq \epsilon_2\}$. Then for all 
$\epsilon,\epsilon_j\leq 1$ we have
\[
\int_{\Delta\setminus\Delta_{\epsilon}^{\alpha}}
\frac{1}{|\zeta_1|\cdots |\zeta_n|}\lesssim
\epsilon^{1/(2|\alpha|)}|\log \epsilon |^{n-1}
\]
and
\[
\int_{\Delta\setminus\Delta_{\epsilon_1,\epsilon_2}^{\alpha, \beta}}
\frac{1}{|\zeta_1|\cdots |\zeta_n|}\lesssim
|(\epsilon_1,\epsilon_2)|^{\omega}, \,\,\, 
2\omega <\mbox{min}\{|\alpha|^{-1}, |\beta|^{-1}\}.
\]
\end{lemma}
\begin{proof}
On the set $\Delta\setminus\Delta_{\epsilon_1,\epsilon_2}^{\alpha, \beta}$,
either $|\zeta^{\alpha}|^2<\epsilon_1$ or $|\zeta^{\beta}|^2<\epsilon_2$ and
so it follows from the first inequality that the integral in the second 
inequality is less then or equal to (a constant times)
\begin{eqnarray*}
\epsilon^{1/(2|\alpha|)}_1|\log \epsilon_1 |^{n-1}+
\epsilon^{1/(2|\beta|)}_2|\log \epsilon_2 |^{n-1} &\lesssim &
\epsilon^{1/(2|\alpha|)-\nu}_1+
\epsilon^{1/(2|\beta|)-\nu}_2 \\
& \lesssim &
|(\epsilon_1,\epsilon_2)|^{\omega_{\nu}},
\end{eqnarray*}
for any $\nu>0$ and 
$\omega_{\nu} \leq \mbox{min}\{|\alpha|^{-1}, |\beta|^{-1}\}/2-\nu$. Hence the 
second inequality follows from the the first one. To prove the first inequality 
we first integrate with respect to the angular variables and then we
make the change of
variables $x_j=\log |\zeta_j|$ to see that the integral in question equals
\begin{equation}\label{eq22}
(4\pi)^n\int_{Q_{\epsilon}}e^{\sum x_j}dx,
\end{equation}
where $Q_{\epsilon}=\{x\in (-\infty,0]^n;\, 2\sum \alpha_jx_j<\log \epsilon\}$.
Since all $x_j\leq 0$ on $Q_{\epsilon}$ 
we have $\exp (\sum x_j)\leq \exp (-|x|)$ here, and choosing 
$R=|\log \epsilon |/(2|\alpha|)$ we see that \eqref{eq22} is less then or equal
to $\int_{\{|x|>R\}}\exp(-|x|)dx$. In polar coordinates this is easily
seen to be of order $\epsilon^{1/(2|\alpha|)}|\log \epsilon |^{n-1}$.
\end{proof}

\begin{lemma}\label{uppsk_lemma1}
Let $\Delta$ be the unit polydisc in $\C^n$ and put 
$\Delta_{\epsilon}^{\alpha}=\{\zeta\in\Delta; |\zeta^{\alpha}|^2\geq \epsilon\}$ and
$\Delta_{\epsilon_1,\epsilon_2}^{\alpha, \beta}=\{\zeta\in\Delta; |\zeta^{\alpha}|^2\geq \epsilon_1,|\zeta^{\beta}|^2\geq \epsilon_2\}$. Then, for $\epsilon,
\epsilon_j \leq 1$, we have
\[
\int_{\Delta_{\epsilon}^{\alpha}}
\frac{\epsilon}{|\zeta^{\alpha}|^2}
\frac{1}{|\zeta_1|\cdots |\zeta_n|} \lesssim
\epsilon^{1/(2|\alpha|)}|\log \epsilon|^{n-1},
\]
\[
\int_{\Delta_{\epsilon_1,\epsilon_2}^{\alpha,\beta}}
\big(\frac{\epsilon_1}{|\zeta^{\alpha}|^2}+\frac{\epsilon_2}{|\zeta^{\beta}|^2}\big)
\frac{1}{|\zeta_1|\cdots |\zeta_n|}
\lesssim |(\epsilon_1,\epsilon_2)|^{\omega}
\]
and
\[
\int_{\Delta_{\epsilon_1,\epsilon_2}^{\alpha,\beta}}
\frac{\epsilon_1 \epsilon_2}{|\zeta^{\alpha}|^2|\zeta^{\beta}|^2}
\frac{1}{|\zeta_1|\cdots |\zeta_n|} 
\lesssim
|(\epsilon_1,\epsilon_2)|^{\omega},
\]
where $2\omega < \mbox{min}\{|\alpha|^{-1}, |\beta|^{-1}\}$.
\end{lemma}
\begin{proof}
The second and third inequality follow from the first one since it implies that the 
integral in the second one is of the size $\epsilon_1^{\tau +1/(2|\alpha|)}+
\epsilon_2^{\tau +1/(2|\beta|)}\lesssim |(\epsilon_1,\epsilon_2)|^{\tau +\omega}$ 
for any $\tau >0$ and 
that the integral in the third is of the size 
$\mbox{min}\{\epsilon_1^{1/(2|\alpha|)}|\log \epsilon_1|^{n-1},
\epsilon_2^{1/(2|\beta|)}|\log \epsilon_2|^{n-1}\}$.
To prove the first inequality we proceed as in the previous lemma and we see that the 
integral in question equals
\begin{eqnarray}\label{eq21}
(4\pi)^n\epsilon\int_{Q_{\epsilon}}
\frac{e^{\sum x_j}}{e^{2\sum \alpha_jx_j}}dx
&=&
(4\pi)^n\epsilon\int_{Q_{\epsilon}\cap \{|x|\leq R\}}
\frac{e^{\sum x_j}}{e^{2\sum \alpha_jx_j}}dx \\
&+&
(4\pi)^n\epsilon\int_{Q_{\epsilon}\cap \{|x|\geq R\}}
\frac{e^{\sum x_j}}{e^{2\sum \alpha_jx_j}}dx,\nonumber
\end{eqnarray}
where $Q_{\epsilon}=\{x\in (-\infty,0]^n;\,2\sum \alpha_jx_j\geq \log \epsilon\}$.
We choose $2R=|\log \epsilon|/|\alpha|$, and then 
$Q_{\epsilon}\cap \{|x|\leq R\}=\{x\in (-\infty,0]^n;\,|x|\leq R\}$. 
If all $x_j\leq 0$ we have $\sum x_j\leq -|x|$ and by the  
Cauchy-Schwarz inequality we also have $-\sum \alpha_jx_j\leq |\alpha| |x|$.
Hence we may estimate the integrand in the second to last integral in 
\eqref{eq21} by $\exp ((2|\alpha|-1)|x|)$. In the last integral we integrate 
where $\epsilon/\exp(2\sum \alpha_jx_j)\leq 1$ and so we see that the right 
hand side of \eqref{eq21} is less then or equal to
\[
(4\pi)^n\epsilon\int_{\{|x|\leq R\}}e^{(2|\alpha|-1)|x|}dx +
(4\pi)^n\int_{\{|x|\geq R\}}e^{-|x|}dx.
\]
By changing to polar coordinates this is seen to be of the size
$\epsilon^{1/(2|\alpha|)}|\log \epsilon |^{n-1}$.
\end{proof}

The proof of the following proposition contains the technical core of this paper.

\begin{prop}\label{tekniskprop}
Assume that $\chi_1,\chi_2 \in C^{\infty}([0,\infty])$ vanish to orders 
$k\geq 0$ and $\ell\geq 0$ at $0$, respectively, and that $\chi_1(\infty)=1$. Then 
for any test form $\varphi \in \D_{n,n}(\C^n)$ we have
\[
\int \frac{1}{\zeta^{k\alpha+\ell\beta}}
\chi_1(\Phi|\zeta^{\alpha}|^2/\epsilon_1)
\chi_2(\Psi|\zeta^{\beta}|^2/\epsilon_2)\varphi \rightarrow
\begin{cases}
\big[\frac{1}{\zeta^{k\alpha+\ell\beta}}\big].\varphi, & \chi_2(\infty)=1 \\
0 , & \chi_2(\infty)=0
\end{cases}
\]
as $\epsilon_1,\epsilon_2\rightarrow 0^+$.
Moreover, as a function of $\epsilon=(\epsilon_1,\epsilon_2)\in [0,\infty)^2$, 
the integral belongs to all $\omega$-Hölder classes with
$2\omega<\mbox{min}\{|\alpha|^{-1},|\beta|^{-1}\}$.
\end{prop}

\begin{remark}
The values of the integral at points $(\epsilon_1,0)$ and $(0,\epsilon_2)$, $\epsilon_j\neq 0$,
are 
\[
\chi_2(\infty)\frac{\chi_1(\Phi|\zeta^{\alpha}|^2/\epsilon_1)}{\zeta^{k\alpha}}
\big[\frac{1}{\zeta^{\ell\beta}}\big].\varphi \,\,\,\mbox{and}\,\,\,
\frac{\chi_2(\Phi|\zeta^{\beta}|^2/\epsilon_2)}{\zeta^{\ell\beta}}
\big[\frac{1}{\zeta^{k\alpha}}\big].\varphi
\]
respectively.
\end{remark}

\begin{remark}
The modulus of continuity can be improved by sharpening the estimates in the
Lemmas \ref{uppsk_lemma2} and \ref{uppsk_lemma1}
but we will not bother about this. This is because the multiindices $\alpha$ and $\beta$
will be implicitly given by Hironaka's theorem and so we can only be sure of the 
existence of some positive Hölder exponent when we prove our main theorems anyway.
\end{remark}

\begin{proof}
We prove H\"older continuity for a path $(\epsilon_1,\epsilon_2)\rightarrow 0$, 
$\epsilon_j\neq 0$. For a general path (inside $[0,\infty)^2$) to an arbitrary point
in $[0,\infty)^2$ one proceeds in a similar way.
Let $K$ be the set of indices $j$ such that $k\alpha_j+\ell\beta_j \geq 2$ and
let $M=M_K^r$ be the operator defined in Lemma \ref{taylorlemma} 
with $r_j=k\alpha_j+\ell\beta_j-2$ for
$j\in K$. Let also $\Delta$ be a polydisc containing the support of $\varphi$.
In this proof we will identify $\varphi$ with its coefficient function with respect
to the volume form in $\C^n$.
We make a preliminary decomposition
\begin{equation}\label{eq06}
\int \frac{1}{\zeta^{k\alpha+\ell\beta}}\chi_1\chi_2\varphi=
\int_{\Delta} \frac{1}{\zeta^{k\alpha+\ell\beta}}\chi_1\chi_2(\varphi-M\varphi)+
\int_{\Delta} \frac{1}{\zeta^{k\alpha+\ell\beta}}\chi_1\chi_2 M\varphi.
\end{equation}
Denote by $\Delta_{\epsilon}$ the set $\{\zeta \in \Delta; |\zeta^{\alpha}|^2 \geq 
\epsilon_1,|\zeta^{\beta}|^2 \geq \epsilon_2\}$. Since $\varphi-M\varphi=
\mathcal{O}(|\zeta^{r+1}|)$, according to Lemma \ref{taylorlemma}, and 
$\chi_1(\infty)=1$ we get
\begin{eqnarray}\label{eq07}
\lefteqn{\Big| \int_{\Delta} 
\frac{1}{\zeta^{k\alpha+\ell\beta}}\chi_1\chi_2(\varphi-M\varphi)-
\chi_2(\infty)
\int_{\Delta} \frac{1}{\zeta^{k\alpha+\ell\beta}}(\varphi-M\varphi)\Big|
}\\
& \lesssim &
\int_{\Delta} \frac{1}{|\zeta_1|\cdots |\zeta_n|}
\big|\chi_1\chi_2-\chi_2(\infty) \big| \nonumber \\
& \lesssim &
\int_{\Delta_{\epsilon}} \frac{1}{|\zeta_1|\cdots |\zeta_n|}
\big|\chi_1\chi_2-\chi_2(\infty) \big|+
\int_{\Delta\setminus\Delta_{\epsilon}} \frac{1}{|\zeta_1|\cdots |\zeta_n|}. \nonumber 
\end{eqnarray} 
It follows from Lemma \ref{uppsk_lemma2} 
that the last integral is of order $|\epsilon|^{\omega}$
as $\epsilon_1,\epsilon_2\rightarrow 0^+$. 
On the other hand, for $\zeta\in \Delta_{\epsilon}$ both 
$|\zeta^{\alpha}|^2/\epsilon_1\geq 1$
and $|\zeta^{\beta}|^2/\epsilon_2\geq 1$ and by Taylor expanding at 
infinity we see that
\begin{eqnarray*}
\chi_1(\Phi|\zeta^{\alpha}|^2/\epsilon_1) 
&=&
\chi_1(\infty)+\frac{\epsilon_1}{|\zeta^{\alpha}|^2}
B_1(\epsilon_1/|\zeta^{\alpha}|^2, \zeta), \\
\chi_2(\Psi|\zeta^{\beta}|^2/\epsilon_2) 
&=&
\chi_2(\infty)+\frac{\epsilon_2}{|\zeta^{\beta}|^2}
B_2(\epsilon_2/|\zeta^{\beta}|^2, \zeta)
\end{eqnarray*}
where $B_1$ and $B_2$ are bounded. Using that $\chi_1(\infty)=1$ we thus get
that $|\chi_1\chi_2-\chi_2(\infty)|$ is of the size
$\epsilon_1/|\zeta^{\alpha}|^2+\epsilon_2/|\zeta^{\beta}|^2$. 
Hence, by Lemma \ref{uppsk_lemma1}
the second to last integral in \eqref{eq07} is also of order
$|\epsilon|^{\omega}$ as $\epsilon_1,\epsilon_2\rightarrow 0^+$.
In view of Lemma \ref{randtermslemma}, we have thus showed that the first
integral on the right hand side of \eqref{eq06} tends to 
$[1/\zeta^{k\alpha+\ell\beta}].\varphi$ if $\chi_2(\infty)=1$ and to zero if
$\chi_2(\infty)=0$ and moreover, belongs to the stated Hölder classes.
We will be done if we can show that
the last integral in \eqref{eq06} is of order $|\epsilon|^{\omega}$.
We know that $M\varphi=\sum_{IJ}\varphi_{IJ}\zeta^I\bar{\zeta}^J$ where each
$\varphi_{IJ}$ is independent of at least one variable and 
$I_j+J_j\leq k\alpha_j+\ell\beta_j-2$ for $j\in K$. Hence, if $\Phi$ and $\Psi$
are constants (or only depend on the modulus of the $\zeta_j$) then the
last integral in \eqref{eq06} is zero for all $\epsilon_1,\epsilon_2 >0$
by anti-symmetry. For the general case, consider one term
\begin{equation}\label{eq08}
\int_{\Delta}\frac{1}{\zeta^{k\alpha+\ell\beta}}\chi_1\chi_2
\varphi_{IJ}\zeta^I\bar{\zeta}^J
\end{equation}
and let $L$ be the set of indices $j\in K$ such that $\zeta_j \mapsto 
\varphi_{IJ}(\zeta)$ is constant. Let also $\mathscr{M}=M_L^{\rho}$ be the 
operator defined in Lemma \ref{taylorlemma} 
with $\rho_j=k\alpha_j+\ell\beta_j-I_j-J_j-2$ for $j\in L$. 
We introduce the independent (real) variables, or ``smoothing parameters'',
$t_1=|\zeta^{\alpha}|^2/\epsilon_1$ and $t_2=|\zeta^{\beta}|^2/\epsilon_2$.
Below, $\mathscr{M}(\chi_1\chi_2)$ denotes the function we obtain by letting 
$\mathscr{M}$ operate on
$\zeta\mapsto \chi_1(t_1\Phi(\zeta))\chi_2(t_2\Psi(\zeta))$ and then 
substituting $|\zeta^{\alpha}|^2/\epsilon_1$ and 
$|\zeta^{\beta}|^2/\epsilon_2$ for $t_1$ and $t_2$ respectively. 
We rewrite the integral \eqref{eq08} as
\begin{eqnarray}\label{eq09}
\int_{\Delta_{\epsilon}}\frac{\varphi_{IJ}\zeta^I\bar{\zeta}^J}{\zeta^{k\alpha+\ell\beta}}
(\chi_1\chi_2-\mathscr{M}(\chi_1\chi_2))&+&
\int_{\Delta\setminus\Delta_{\epsilon}}
\frac{\varphi_{IJ}\zeta^I\bar{\zeta}^J}{\zeta^{k\alpha+\ell\beta}}
(\chi_1\chi_2-\mathscr{M}(\chi_1\chi_2))\nonumber\\
&+&
\int_{\Delta}\frac{\varphi_{IJ}\zeta^I\bar{\zeta}^J}{\zeta^{k\alpha+\ell\beta}}
\mathscr{M}(\chi_1\chi_2).
\end{eqnarray}
Now, $\mathscr{M}(\chi_1\chi_2)$ is a sum of 
terms which, at least for some $j\in L$, are monomials in $\zeta_j$ and 
$\bar{\zeta}_j$ times coefficient functions depending on $|\zeta_j|$ and 
the other variables. The degrees of these monomials do not exceed 
$\rho_j=k\alpha_j+\ell\beta_j-I_j-J_j-2$ and since 
$\zeta_j\mapsto \varphi_{IJ}(\zeta)$ is constant 
for $j\in L$ we see, by counting exponents, that the last integral in 
\eqref{eq09} vanishes by anti-symmetry for all $\epsilon_1,\epsilon_2>0$.
By Lemma \ref{trollerilemma} we have
\begin{equation}\label{eq10}
\chi_1(t_1\Phi)\chi_2(t_2\Psi)-
\mathscr{M}(\chi_1(t_1\Phi)\chi_2(t_2\Psi))=
|\zeta^{\rho+1}|B(t_1,t_2,\zeta),
\end{equation}
where $B$ is bounded on $(0,\infty)^2\times \Delta$. We note also that 
by Lemma \ref{taylorlemma},
$\varphi_{IJ}(\zeta)=\mathcal{O}(\prod_{j\in L\setminus K}|\zeta_j|^{r_j+1})$.
From \eqref{eq10} we thus see that
the modulus of the second integral 
in \eqref{eq09} can be estimated by
\[
C\int_{\Delta\setminus\Delta_{\epsilon}}\frac{1}{|\zeta_1|\cdots|\zeta_n|},
\]
which is of order $|\epsilon|^{\omega}$ by Lemma \ref{uppsk_lemma2}. 
It remains to consider the 
first integral in \eqref{eq09}. On the set $\Delta_{\epsilon}$ we have that
$\Phi|\zeta^{\alpha}|^2/\epsilon_1$ and $\Psi|\zeta^{\beta}|^2/\epsilon_2$
are larger then some positive constant and so by multiplying the Taylor expansions of
the functions $t_1\mapsto \chi_1(t_1\Phi)$ and $t_2\mapsto \chi_2(t_2\Psi)$ at infinity
we get
\begin{eqnarray*}
\chi_1(\Phi|\zeta^{\alpha}|^2/\epsilon_1)
\chi_2(\Psi|\zeta^{\beta}|^2/\epsilon_2)&=&
\chi_2(\infty)+
\frac{\epsilon_2}{|\zeta^{\beta}|^2}
\tilde{\chi}_2(|\zeta^{\beta}|^2/\epsilon_2,\zeta)\\
&+&
\chi_2(\infty)\frac{\epsilon_1}{|\zeta^{\alpha}|^2}
\tilde{\chi}_1(|\zeta^{\alpha}|^2/\epsilon_1,\zeta)\\
&+&
\frac{\epsilon_1\epsilon_2}{|\zeta^{\alpha}|^2|\zeta^{\beta}|^2}
\tilde{\chi}_1(|\zeta^{\alpha}|^2/\epsilon_1,\zeta)
\tilde{\chi}_2(|\zeta^{\beta}|^2/\epsilon_2,\zeta)
\end{eqnarray*}
where $\tilde{\chi}_j$ are smooth on $[1,\infty]\times \Delta$. Now since 
$|\zeta^{\alpha}|^2/\epsilon_1=t_1$ and 
$|\zeta^{\beta}|^2/\epsilon_2=t_2$ are independent variables we conclude that
\begin{eqnarray*}
\chi_1\chi_2-\mathscr{M}(\chi_1\chi_2) &=& 
\frac{\epsilon_2}{|\zeta^{\beta}|^2}(\tilde{\chi}_2-\mathscr{M}\tilde{\chi}_2)
+
\frac{\epsilon_1}{|\zeta^{\alpha}|^2}\chi_2(\infty)
(\tilde{\chi}_1-\mathscr{M}\tilde{\chi}_1)
\\
&+&
\frac{\epsilon_1\epsilon_2}{|\zeta^{\alpha}|^2|\zeta^{\beta}|^2}
(\tilde{\chi}_1\tilde{\chi}_2-\mathscr{M}(\tilde{\chi}_1\tilde{\chi}_2))
\end{eqnarray*}
for $\zeta \in \Delta_{\epsilon}$. By Lemmas \ref{taylorlemma} and \ref{uppsk_lemma1}
we see that the first integral in \eqref{eq09} also is of order 
$|\epsilon|^{\omega}$ as $\epsilon_1,\epsilon_2\rightarrow 0^+$ 
and the proof is complete.
\end{proof}

\begin{remark}\label{SEPrem}
Let us assume that the function $\Phi$ is identically $1$ in the previous proposition.
Then, instead of adding and subtracting $\mathscr{M}(\chi_1\chi_2)$ in \eqref{eq09},
it is enough to add and subtract $\chi_1\mathscr{M}(\chi_2)$.
This suggests that one can relax the smoothness assumption on $\chi_1$. 
It is actually possible to take $\chi_1$ to be the characteristic function of $[1,\infty]$. 
If we define the value of the integral in Proposition \ref{tekniskprop} at a point 
$(\epsilon_1,0)$ to be 
\begin{equation}\label{schwartzprod}
\int_{\Delta}\frac{1}{\zeta^{k\alpha+\ell\beta}}\chi_1(|\zeta^{\alpha}|^2/\epsilon_1)
(\varphi-M\varphi),
\end{equation}
where $\Delta$ and $M$ are as in the proof above, then the conclusions of Proposition
\ref{tekniskprop} hold for this choice of $\chi_1$. Only minor changes in the proof are needed
to see this.
One can also check that \eqref{schwartzprod} is a way of computing 
\[
\chi_1(|\zeta^{\alpha}|^2/\epsilon_1)
\big[\frac{1}{\zeta^{k\alpha+\ell\beta}}\big].\varphi.
\]
The product $\chi_1(|\zeta^{\alpha}|^2/\epsilon_1)[1/\zeta^{k\alpha+\ell\beta}]$ is well
defined because the wave front sets of the two currents 
behave in the right way, at least for almost all $\epsilon_1$, see \cite{jeb2}.
\end{remark}

We make another useful observation.
Since the function $\tilde{\chi}(s)= \chi(1/s)$ is smooth at 
zero and 
$\tilde{\chi}'(s):=-\frac{1}{s^2}\chi'(1/s)$,
it follows that $s\mapsto \chi'(1/s)/s$ is smooth at zero and vanishes for $s=0$. 
Hence, $t\mapsto \chi'(t)t$ is smooth on $[0,\infty]$, vanishes
to the same order at zero as $\chi$, and maps $\infty$ to $0$.
From Proposition \ref{tekniskprop} we thus see that we have

\begin{cor}\label{korret}
Assume that $\chi_1,\chi_2\in C^{\infty}([0,\infty])$ vanish to orders 
$k$ and $\ell$ at zero respectively, and satisfy $\chi_j(\infty)=1$. For any 
smooth and strictly positive functions $\Phi$ and $\Psi$ on $\C^n$ and any
test form $\varphi\in \D_{n,n}(\C^n)$ we have
\begin{equation}
\lim_{\epsilon_1,\epsilon_2\rightarrow 0^+}\int
\frac{1}{\zeta^{k\alpha+\ell\beta}}
\chi_1(\Phi|\zeta^{\alpha}|^2/\epsilon_1)
\chi'_2(\Psi|\zeta^{\beta}|^2/\epsilon_2)
\frac{|\zeta^{\beta}|^2}{\epsilon_2}\varphi
=0,
\end{equation}
and moreover, as a function of 
$\epsilon=(\epsilon_1,\epsilon_2)\in [0,\infty)^2$, the integral belongs to all
$\omega$-Hölder classes with $2\omega<\mbox{min}\{|\alpha|^{-1},|\beta|^{-1}\}$.
\end{cor}

\section{Regularizations of products of Cauchy-Fantappiè-Leray type currents}\label{prodsektion}
We are now in a position to prove our main results. We start with a regularization
of the product $U^f\wedge U^g$. Recall that if $f$ is function then 
$U^f=[1/f]$ times some basis element.
 
\begin{thm}\label{korv}
Let $f$ and $g$ be holomorphic sections (locally non-trivial) of the 
holomorphic $m_j$-bundles $E_j^*\rightarrow X$, $j=1,2$, respectively. Let 
$\chi_1,\chi_2\in C^{\infty}([0,\infty])$ be any functions vanishing to orders
$m_1$ and $m_2$ at zero respectively, and satisfying $\chi_j(\infty)=1$. Then, 
for any test form $\varphi$ we have
\[
\int
\chi_1(|f|^2/\epsilon_1)u^f\wedge\chi_2(|g|^2/\epsilon_2)u^g\wedge \varphi
\rightarrow
U^f\wedge U^g.\varphi, 
\]
as $\epsilon_1,\epsilon_2\rightarrow 0^+$. Moreover, as a function of 
$\epsilon=(\epsilon_1,\epsilon_2)\in [0,\infty)^2$
the integral on the left hand side belongs to some Hölder class independently of $\varphi$.
\end{thm}

\begin{proof}
Recall that $U^f\wedge U^g.\varphi$ is defined as the value at zero of 
the meromorphic function
\[
\lambda \mapsto 
\int |f|^{2\lambda}u^f\wedge |g|^{2\lambda}u^g\wedge \varphi.
\]
Assuming only that $\chi_1$ and $\chi_2$ vanish to orders $k\leq m_1$ and 
$\ell\leq m_2$ at zero respectively we will show that 
\begin{equation}\label{eql}
\int \chi_1u^f_{k,k-1}\wedge \chi_2 u^g_{\ell,\ell-1}\wedge \varphi
\rightarrow 
\int |f|^{2\lambda}u^f_{k,k-1}\wedge |g|^{2\lambda} u^g_{\ell,\ell-1}\wedge \varphi
\Big|_{\lambda=0}
\end{equation}
and that the left hand side belongs to some Hölder class. This will clearly 
imply the theorem. We may assume that $\varphi$ has arbitrarily small support 
after a partition of unity. If $\varphi$ has support outside 
$f^{-1}(0)\cup g^{-1}(0)$ it is easy to check that \eqref{eql} holds and hence 
we can restrict to the case that $\varphi$ has support in a small neighborhood 
$\mathcal{U}$ of a point $p\in f^{-1}(0)\cup g^{-1}(0)$.
We may also assume that $\mathcal{U}$ is contained in a coordinate neighborhood
and that all bundles are trivial over $\mathcal{U}$. We let 
$(f_1,\ldots,f_{m_1})$ and $(g_1,\ldots,g_{m_2})$ denote the components of
$f$ and $g$ respectively, with respect to some holomorphic frames. It follows 
from Hironaka's theorem, possibly after another localization, that there is 
an $n$-dimensional complex manifold $\tilde{\mathcal{U}}$ and a proper 
holomorphic map $\Pi\colon \tilde{\mathcal{U}}\rightarrow \mathcal{U}$ such
that $\Pi$ is biholomorphic outside the nullset 
$\Pi^*\{f_1\cdots f_{m_1}\cdot g_1\cdots g_{m_2}\}$ and that this hypersurface
has normal crossings in $\tilde{\mathcal{U}}$. Hence we can cover 
$\tilde{\mathcal{U}}$ by local charts, each centered at the origin,
such that $\Pi^*f_j$ and $\Pi^*g_j$ are monomials times non-vanishing functions.
The support of $\Pi^*\varphi$ is compact because $\Pi$ is proper and hence, 
we can cover the support of $\Pi^*\varphi$ by finitely many of these charts. 
We let $\rho_k$ be a partition of unity on 
$\mbox{supp}(\Pi^*\varphi)$ subordinate to this cover. 
Now, following \cite{PTY} and \cite{BGVY}, given monomials $\mu_1\ldots,\mu_{\nu}$, one
can construct an $n$-dimensional toric manifold $\mathcal{X}$ and a proper 
holomorphic map $\tilde{\Pi}\colon \mathcal{X}\rightarrow \C^n_t$ which is 
monoidal when expressed in local coordinates in each chart. Moreover, 
$\tilde{\Pi}$ is biholomorphic outside $\tilde{\Pi}^*\{t_1\cdots t_n=0\}$ and in 
each chart one of the monomials 
$\tilde{\Pi}^*\mu_1,\ldots,\tilde{\Pi}^*\mu_{\nu}$ divides all the others.
By repeating this process, if necessary, and localizing with partitions of unity
at each step, we may actually assume that $f_j=\mu_{f,j}\tilde{f}_j$ and 
$g_j=\mu_{g,j}\tilde{g}_j$ where $\tilde{f}_j$ and $\tilde{g}_j$ are 
non-vanishing and $\mu_{f,j}$ and $\mu_{g,j}$ are monomials with the property
that $\mu_{f,\nu_1}$ divides all $\mu_{f,j}$ and $\mu_{g,\nu_2}$ divides all
$\mu_{g,j}$ for some indices $\nu_1$ and $\nu_2$. Denote $\mu_{f,\nu_1}$ by
$\zeta^{\alpha}$ and $\mu_{g,\nu_2}$ by $\zeta^{\beta}$. It follows that 
$|f|^2=|\zeta^{\alpha}|^2\Phi$ and $|g|^2=|\zeta^{\beta}|^2\Psi$ where 
$\Phi$ and $\Psi$ are strictly positive functions. Moreover, 
$s_f=\bar{\zeta}^{\alpha}\tilde{s}_f$ and 
\[
u^f_{k,k-1}=\frac{s_f\wedge (\debar s_f)^{k-1}}{|f|^{2k}}=
\frac{1}{\zeta^{k\alpha}}
\frac{\tilde{s}_f\wedge (\debar \tilde{s}_f)^{k-1}}{\Phi^k}=
\frac{1}{\zeta^{k\alpha}}\tilde{u}^f_{k,k-1}
\]
where $\tilde{u}^f_{k,k-1}$ is a smooth form and similarly for $u^g_{\ell,\ell-1}$. 
In order to prove \eqref{eql} it thus suffices to prove
\begin{eqnarray}\label{eqi}
\lefteqn{\int \frac{\chi_1(\Phi |\zeta^{\alpha}|^2/\epsilon_1)}{\zeta^{k\alpha}}
\tilde{u}^f_{k,k-1}
\wedge \frac{\chi_2(\Psi |\zeta^{\beta}|^2/\epsilon_2)}{\zeta^{\ell\beta}} 
\tilde{u}^g_{\ell,\ell-1}\wedge 
\tilde{\varphi}}\\
& &
\rightarrow
\int 
\frac{|\zeta^{\alpha}|^{2\lambda}}{\zeta^{k\alpha}}
\Phi^{\lambda}\tilde{u}^f_{k,k-1}\wedge 
\frac{|\zeta^{\beta}|^{2\lambda}}{\zeta^{\ell\beta}}
\Psi^{\lambda}\tilde{u}^g_{\ell,\ell-1}\wedge 
\tilde{\varphi}
\Big|_{\lambda=0} \nonumber
\end{eqnarray}
where $\tilde{\varphi}=\rho_{k_j}\Pi_j^*\cdots\rho_{k_1}\Pi_1^*\varphi$ and that 
the integral on the left hand side belongs to some Hölder class.
But by 
Proposition \ref{tekniskprop} it does belong to some Hölder class and tends 
to $[1/\zeta^{k\alpha+\ell\beta}].\tilde{u}^f_{k,k-1}\wedge \tilde{u}^g_{\ell,\ell-1}\wedge\tilde{\varphi}$. One can verify that this indeed is equal to the right 
hand side of \eqref{eqi} by integrations by parts as in e.g.\ \cite{A1}.
\end{proof}

\begin{remark}
This theorem can actually be generalized to any number of factors $U^f$. One first checks
that the analogue of Proposition \ref{tekniskprop} holds for any number of 
functions $\chi_j$ and then reduces to this case just as in the proof above. 
In particular, if $f_j$, $j=1,\ldots,p$, are holomorphic functions and 
$\chi_j$ vanish at $0$, we have
\[
\int \frac{\chi_1(|f_1|^2/\epsilon_1)}{f_1}\cdots 
\frac{\chi_p(|f_p|^2/\epsilon_p)}{f_p}\,\varphi \rightarrow
\big[\frac{1}{f_1}\cdots \frac{1}{f_p}\big].\varphi
\]
unrestrictedly as all $\epsilon_j\rightarrow 0^+$.
However, we focus on the two factor case since
we do not know how to handle more than two residue factors.
\end{remark}

To prove our regularization results for the currents $U^f\wedge R^g$ and $R^f\wedge R^g$ 
we have to structure the information obtained from an application of Hironaka's
theorem more carefully 
and then use Proposition \ref{tekniskprop} and Corollary \ref{korret} in the right 
way. The technical part of this is contained in the following proposition.

\begin{prop}\label{palle}
Assume that $\chi_1,\chi_2\in C^{\infty}([0,\infty])$ vanish to orders 
$k$ and $\ell$ at zero, respectively, and satisfy $\chi_j(\infty)=1$. Let
$\alpha'$, $\alpha''$, $\beta'$ and $\beta''$ be 
multiindices such that $\alpha'$, 
$\alpha''$ and $\beta'$ have pairwise disjoint supports, and $\alpha''_j=0$ if
and only if $\beta''_j=0$. Assume also that $\varphi\in \D_{n,n-1}(\C^n)$
has the property that 
$d\bar{\zeta}_j/\bar{\zeta}_j \wedge \varphi\in \D_{n,n}(\C^n)$ for 
all $j$ such that $\alpha''_j\neq 0$. Then for any 
smooth and strictly positive functions $\Phi$ and $\Psi$ on $\C^n$ we have
\[
\lim_{\epsilon_1,\epsilon_2 \rightarrow 0^+}
\int \frac{1}{\mu_1^k\mu_2^\ell}
\chi_1(\Phi|\mu_1|^2/\epsilon_1)
\debar\chi_2(\Psi|\mu_2|^2/\epsilon_2)\wedge \varphi =
\Big[\frac{1}{\mu_1^k\zeta^{\ell\beta''}}\Big]\otimes \debar
\Big[\frac{1}{\zeta^{\ell\beta'}}\Big].\varphi,
\]
where $\mu_1=\zeta^{\alpha'+\alpha''}$ and $\mu_2=\zeta^{\beta'+\beta''}$.
Moreover, as a function of 
$\epsilon=(\epsilon_1,\epsilon_2)\in [0,\infty)^2$, the integral belongs to all
$\omega$-Hölder classes with $2\omega<\mbox{min}\{|\alpha'+\alpha''|^{-1},
|\beta'+\beta''|^{-1}\}$.
\end{prop}

\begin{remark}
Note that the hypotheses on the multiindices imply that a factor $\zeta_j$
divides both the monomials $\mu_1$ and $\mu_2$ if and only if $\alpha''_j\neq 0$
(or equivalently $\beta''_j\neq 0$). In particular, the tensor product of 
the currents is well defined.
\end{remark}
\begin{remark}\label{rem1}
We may let $k$ or $\ell$ or both of them be equal to zero and the conclusions of 
the proposition still hold. In case $\ell=0$ one should interpret 
$\debar [1/\zeta^{\ell\beta'}]$ as zero. 
\end{remark}

\begin{proof}
Let $K$, $L$ and $K^c$ be the set of indices $j$ such that $\beta'_j\neq 0$,
$\beta''_j\neq 0$ and $\beta'_j=0$ respectively. Clearly $L\subseteq K^c$.
We write $\debar=\debar_K+\debar_{K^c}$ and integrate by parts with respect to
$\debar_K$ to see that
\begin{eqnarray}\label{eqa}
\lefteqn{\int
\frac{1}{\mu_1^k\mu_2^\ell}\chi_1(\debar_K+\debar_{K^c})\chi_2 \wedge \varphi =}\\
&- &
\int \frac{1}{\mu_1^k\mu_2^\ell}\chi'_1\frac{|\mu_1|^2}{\epsilon_1}\chi_2
\debar_K \Phi \wedge \varphi
-\int \frac{1}{\mu_1^k\mu_2^\ell}\chi_1\chi_2 \debar_K \varphi \nonumber \\
&+&
\int \frac{1}{\mu_1^k\mu_2^\ell}\chi_1\chi'_2\frac{|\mu_2|^2}{\epsilon_2}
(\Psi\sum_{j\in L}\beta''_j\frac{d\bar{\zeta}_j}{\bar{\zeta}_j}+
\debar_{K^c}\Psi)\wedge \varphi.\nonumber
\end{eqnarray}
Note that $\debar_K$ does not fall on $|\mu_1|^2$ because of the hypotheses on 
the multiindices. By assumption, 
$d\bar{\zeta}_j/\bar{\zeta}_j \wedge \varphi\in \D_{n,n}(\C^n)$ for $j\in L$ and
so the first and the last integral on the right hand side of \eqref{eqa}
tend to zero and has the right modulus of continuity by Corollary \ref{korret}.
The second to last integral in \eqref{eqa} tends to 
$-[1/(\mu_1^k\mu_2^\ell)].\debar_K \varphi=[1/(\mu_1^k\zeta^{\ell\beta''})]\otimes 
\debar [1/\zeta^{\ell\beta'}].\varphi$ and has the right modulus of continuity
by Proposition \ref{tekniskprop}.
\end{proof}

\begin{thm}\label{satsen}
Let $f$ and $g$ be holomorphic sections (locally non-trivial) of the 
holomorphic $m_j$-bundles $E_j^*\rightarrow X$, $j=1,2$, respectively. 
Assume that the section $f\oplus g$
of $E_1^*\oplus E_2^*\rightarrow X$ defines a complete intersection. Let 
$\chi_1,\chi_2\in C^{\infty}([0,\infty])$ be any functions vanishing to orders
$m_1$ and $m_2$ at zero respectively, and satisfying $\chi_j(\infty)=1$. Then, 
for any test form $\varphi$ we have
\begin{equation}\label{eqb}
\int 
\chi_1(|f|^2/\epsilon_1)u^f\wedge \debar \chi_2(|g|^2/\epsilon_2)\wedge u^g
\wedge \varphi \rightarrow
U^f\wedge R^g.\varphi
\end{equation}
as $\epsilon_1,\epsilon_2\rightarrow 0^+$. Moreover, as a function of 
$\epsilon=(\epsilon_1,\epsilon_2)\in [0,\infty)^2$ the 
integral on the left hand side belongs to some Hölder class independently of $\varphi$.
\end{thm}

\begin{proof}
We will assume that $\chi_1$ and $\chi_2$
only vanish to orders $k\leq m_1$ and $\ell\leq m_2$ respectively and show that
\begin{equation}\label{eqc}
\int
\chi_1u^f_{k,k-1}\wedge \debar \chi_2\wedge u^g_{\ell,\ell-1}\wedge \varphi \rightarrow
\int |f|^{2\lambda}u^f_{k,k-1}\wedge\debar |g|^{2\lambda}\wedge 
u^g_{\ell,\ell-1}\wedge \varphi\Big|_{\lambda=0}.
\end{equation}
By arguing as in the proof of Theorem \ref{korv} we may assume that 
$|f|^2=|\zeta^{\alpha}|^2\Phi$ and $|g|^2=|\zeta^{\beta}|^2\Psi$ where 
$\Phi$ and $\Psi$ are strictly positive functions and moreover, that 
$u^f_{k,k-1}=\tilde{u}^f_{k,k-1}/\zeta^{k\alpha}$ for a smooth form
$\tilde{u}^f_{k,k-1}$ and similarly for $u^g_{\ell,\ell-1}$. What we have to prove 
is thus
\begin{eqnarray}\label{eqj}
\lefteqn{\int \frac{\chi_1(\Phi |\zeta^{\alpha}|^2/\epsilon_1)}{\zeta^{k\alpha}}
\tilde{u}^f_{k,k-1}
\wedge \frac{\debar\chi_2(\Psi |\zeta^{\beta}|^2/\epsilon_2)}{\zeta^{\ell\beta}} 
\tilde{u}^g_{\ell,\ell-1}\wedge 
\tilde{\varphi}}\\
& &
\rightarrow
\int 
\frac{|\zeta^{\alpha}|^{2\lambda}}{\zeta^{k\alpha}}
\Phi^{\lambda}\tilde{u}^f_{k,k-1}\wedge 
\frac{\debar(|\zeta^{\beta}|^{2\lambda}\Psi^{\lambda})}{\zeta^{\ell\beta}}
\tilde{u}^g_{\ell,\ell-1}\wedge 
\tilde{\varphi}
\Big|_{\lambda=0} \nonumber
\end{eqnarray}
where $\tilde{\varphi}=\rho_{k_j}\Pi_j^*\cdots\rho_{k_1}\Pi_1^*\varphi$.
After the resolutions of singularities we can in general no longer say that 
the pull-back of $f\oplus g$ defines a complete intersection. On the other hand 
we claim that if $\zeta_j$ divides both $\zeta^{\alpha}$ and $\zeta^{\beta}$
then $d\bar{\zeta}_j/\bar{\zeta}_j\wedge\tilde{\varphi}$ is smooth.
In fact, let $z$ be local coordinates on our original manifold. In order that
the integrals in \eqref{eqc} should be non-zero, $\varphi$ has to have degree
$n-k-\ell+1$ in $d\bar{z}$ and so we can assume that 
\[
\varphi=\sum_{\#J=n-k-\ell+1}\varphi_J\wedge d\bar{z}_J.
\]
Since the variety $V_{f\oplus g}=f^{-1}(0)\cap g^{-1}(0)$ has dimension
$n-m_1-m_2<n-k-\ell+1$ we see that $d\bar{z}_J$ vanishes on $V_{f\oplus g}$. 
The pull-back of $d\bar{z}_J$ through all the resolutions $\Pi_j$ can be written
$\sum_I C_I(\zeta)d\bar{\zeta}_I$ and it must vanish on the pull-back of 
$V_{f\oplus g}$. In particular it has to vanish on $\{\zeta_j=0\}$ if 
$\zeta_j$ divides both $\zeta^{\alpha}$ and $\zeta^{\beta}$.
If $d\bar{\zeta}_j$ does not occur in $d\bar{\zeta}_I$ 
it must be that the coefficient
function $C_I(\zeta)$ vanishes on $\{\zeta_j=0\}$. But these functions are 
anti-holomorphic and so $\bar{\zeta}_j$ must divide $C_I(\zeta)$. The claim is 
established. We now write $\zeta^{\alpha}=\zeta^{\alpha'+\alpha''}$ and
$\zeta^{\beta}=\zeta^{\beta'+\beta''}$ where $\alpha'$, $\alpha''$ and $\beta'$
have pairwise disjoint supports and $\alpha''=0$ if and only if $\beta''=0$. 
Thus, $\zeta_j$ divides both $\zeta^{\alpha}$ and $\zeta^{\beta}$
if and only if $\alpha''_j\neq 0$, or equivalently, $\beta''_j\neq 0$.
According to Proposition \ref{palle} the left hand side of \eqref{eqj} 
belongs to some Hölder class and tends to
\[
-\Big[\frac{1}{\zeta^{k\alpha+\ell\beta''}}\Big]\otimes
\debar\Big[\frac{1}{\zeta^{\ell\beta'}}\Big].
\tilde{u}^f_{k,k-1}\wedge\tilde{u}^g_{\ell,\ell-1}\wedge \tilde{\varphi}.
\]
One can compute the right hand side of \eqref{eqj} by integrations by parts 
as in e.g.\ \cite{A1} to see that it equals the same thing.
\end{proof}

\begin{remark}\label{morot}
The form $\debar \chi_2(|g|^2/\epsilon_2)\wedge u^g$ is actually smooth even if $\chi_2$ only
vanishes to order $m_2$ at $0$. The only possible problem is with the top degree term
$\debar \chi_2(|g|^2/\epsilon_2)\wedge u_{m_2,m_2-1}^g$. But we have
\begin{eqnarray*}
C^{\infty}(X)\ni \debar (\chi_2(|g|^2/\epsilon_2)u_{m_2,m_2-1}^g)
&=&
\debar \chi_2(|g|^2/\epsilon_2)\wedge u_{m_2,m_2-1}^g\\
&+&
\chi_2(|g|^2/\epsilon_2)\debar u_{m_2,m_2-1}^g,
\end{eqnarray*}
and since $u_{m_2,m_2-1}^g$ is $\debar$-closed (outside $V_g$) it follows that  
$\debar \chi_2(|g|^2/\epsilon_2)\wedge u_{m_2,m_2-1}^g$ is smooth as well.
\end{remark}

\begin{cor}\label{sven}
With the same hypotheses as in Theorem \ref{satsen} we have

\[
\int
\debar \chi_1(|f|^2/\epsilon_1)\wedge u^f\wedge
\debar \chi_2(|g|^2/\epsilon_2)\wedge u^g\wedge\varphi 
\rightarrow R^f\wedge R^g .\varphi,
\]
\begin{equation}\label{eqbb}
\int 
\debar \chi_1(|f|^2/\epsilon_1)\wedge u^f
\chi_2(|g|^2/\epsilon_2)\wedge \varphi \rightarrow
R^f.\varphi,
\end{equation}
and
\begin{equation}\label{eqbbb}
\int 
\chi_1(|f|^2/\epsilon_1)\wedge u^f\wedge\debar\chi_2(|g|^2/\epsilon_2)
\wedge \varphi \rightarrow
0
\end{equation}
as $\epsilon_1,\epsilon_2\rightarrow 0^+$, and as functions of 
$\epsilon=(\epsilon_1,\epsilon_2)\in [0,\infty)^2$ the 
integrals on the left hand sides belong to some Hölder classes independently of $\varphi$.
\end{cor} 
\begin{proof}
We have the following equality of smooth forms:
\begin{eqnarray}\label{eqna}
\nabla(\debar\chi_1 \wedge u^f\wedge \chi_2 u^g)&=&
-\debar \chi_1\wedge \chi_2 u^g 
-\debar\chi_1\wedge u^f\wedge \debar\chi_2\wedge u^g\\
&+&\debar\chi_1 \wedge u^f\chi_2. \nonumber
\end{eqnarray}
The computation rules established in \cite{W}, and
Theorem \ref{satsen} now imply that, for any test form $\varphi$ (of complementary
total degree), we have
\begin{eqnarray*}
R^f.\varphi-R^f\wedge R^g.\varphi &=&
\nabla (R^f\wedge U^g).\varphi = -R^f\wedge U^g.\nabla\varphi\\ 
&=&
\lim -\int\debar\chi_1 \wedge u^f\wedge \chi_2 u^g\wedge \nabla\varphi\\
&=& \lim \int \nabla(\debar\chi_1 \wedge u^f\wedge \chi_2 u^g)\wedge \varphi.
\end{eqnarray*}
The integral on the second row is Hölder continuous by Theorem \ref{satsen}
and so, also the integral on the third row is. By choosing $\varphi$ of appropriate
bidegrees the corollary now follows from \eqref{eqna}.
\end{proof}

The statements \eqref{eqbb} and \eqref{eqbbb} actually hold with no assumptions 
on the behavior of $\chi_2$ at zero. This can be seen by using that we know
this when $\chi_2\equiv 1$ by Corollary \ref{spenat}, and when $\chi_2$ 
vanishes to high enough order by the previous corollary.

Assume that $f$ defines a complete intersection and
pick a holomorphic function $g$ such that $f\oplus g$ also defines a complete intersection
and such that $g$ is zero on the singular part of $V_f$. 
After resolving singularities in the 
proof of Theorem \ref{satsen} we can find coordinates such that $g$ is a monomial times 
a non-vanishing holomorphic function $\tilde{g}$. 
But $\tilde{g}$ can be incorporated in some coordinate
and we can therefore assume that $\tilde{g}\equiv 1$. 
Repeating the proof of Theorem \ref{satsen} and using Remark \ref{SEPrem} one shows that 
\eqref{eqbb} holds for $\chi_2$ equal to the characteristic function of $[1,\infty]$.
Then, if we first let $\epsilon_1$ tend to zero, keeping $\epsilon_2$ fixed, and after that 
let $\epsilon_2$ tend to zero we get that 
\[
\lim_{\epsilon_2\rightarrow 0^+}\chi_2(|g|^2/\epsilon_2)R^f=R^f.
\]
We remark that the product $\chi_2(|g|^2/\epsilon_2)R^f$ is well defined since the 
wave front sets of 
$\chi_2(|g|^2/\epsilon_2)$ and $R^f$ behave properly, see e.g.\ \cite{jeb2}. 
Since $\chi_2(|g|^2/\epsilon_2)$ equals the characteristic function of $\{|g|^2>\epsilon_2\}$
we have 

\begin{cor}\label{SEP}
If $f$ defines a complete intersection then the Cauchy-Fantappiè-Leray current $R^f$
has the standard extension property.
\end{cor}

This is a well known result and follows from the fact that $R^f$ equals the 
Coleff-Herrera current in the sense of \eqref{BM=CH}. 
It is even true that $\chi_{\rho g}(\epsilon)R^f\rightarrow R^f$,
$\epsilon\rightarrow 0^+$ where $\rho$ is a positive smooth function and
$\chi_{\rho g}(\epsilon)$ is the characteristic function of $\{|\rho g|>\epsilon\}$. In fact,
via Hironaka and toric resolutions one reduces to the case of one function and then
one can proceed as in \cite{jeb2}.

We know from \cite{W} that if $f\oplus g$ defines a complete intersection then 
$R^f\wedge R^g$ consists of one term of top degree. Hence, it is only the top degree
term of $\debar \chi_1\wedge u^f\wedge \debar \chi_2\wedge u^g$ which gives a 
contribution in the limit.
With the natural choices $\chi_1(t)=t^{m_1}/(t+1)^{m_1}$ and 
$\chi_2(t)=t^{m_2}/(t+1)^{m_2}$, Corollary \ref{sven} and Remark \ref{morot} thus give

\begin{cor}\label{leif}
Let $f$ and $g$ be holomorphic sections (locally non-trivial) of the 
holomorphic $m_j$-bundles $E_j^*\rightarrow X$, $j=1,2$, respectively. 
Assume that the section $f\oplus g$
of $E_1^*\oplus E_2^*\rightarrow X$ defines a complete intersection. Then, 
for any test form $\varphi$ we have
\[
\int
\debar \frac{s_f\wedge(\debar s_f)^{m_1-1}}{(|f|^2+\epsilon_1)^{m_1}}\wedge
\debar \frac{s_g\wedge(\debar s_g)^{m_2-1}}{(|g|^2+\epsilon_2)^{m_2}}\wedge
\varphi \rightarrow
R^f\wedge R^g.\varphi
\]
as $\epsilon_1,\epsilon_2\rightarrow 0^+$, and the integral to the left belongs
to some Hölder class independently of $\varphi$.
\end{cor}

For sections $f$ and $g$ of the trivial line bundle we get the result announced
in \cite{HS}.

\begin{cor}\label{CR}
Let $f$ and $g$ be holomorphic functions defining a complete intersection.
Then for any test form $\varphi$ we have
\[
\int
\debar \frac{\bar{f}}{|f|^2+\epsilon_1}\wedge
\debar \frac{\bar{g}}{|g|^2+\epsilon_2}\wedge \varphi \rightarrow
\Big[\debar\frac{1}{f}\wedge\debar\frac{1}{g}\Big].\varphi
\]
as $\epsilon_1,\epsilon_2\rightarrow 0^+$, and the integral to the left 
belongs to some Hölder class independently of $\varphi$.
\end{cor}

\begin{proof}
We consider $f$ and $g$ as sections of (different copies of)
the trivial line bundle 
$X\times\C\rightarrow X$ with the standard metric. Then, suppressing the natural
global frame elements, we have $s_f=\bar{f}$ and $s_g=\bar{g}$. By Corollary
\ref{leif} we are done since $R^f\wedge R^g$ is the Coleff-Herrera current.
\end{proof}

So far, in this section, we have used one function $\chi$ to regularize all terms of
$u^f$. One could try to take different $\chi$:s for different terms. We recall the 
natural choices $t^{k}/(t+1)^{k}$ from
Corollary \ref{sture} and we let $u^f_{\epsilon}=s_f/(\nabla s_f+\epsilon)=
\sum s_f\wedge (\debar s_f)^{k-1}/(|f|^2+\epsilon)^k$.
The next theorem says that, in the complete intersection case,  
the product of two such regularized currents goes unrestrictedly to the product, 
in the sense of \cite{W}, of the currents.

\begin{thm}\label{kent}
Let $f$ and $g$ be holomorphic sections (locally non-trivial) of the 
holomorphic $m_j$-bundles $E_j^*\rightarrow X$, $j=1,2$, respectively. 
Assume that the section $f\oplus g$
of $E_1^*\oplus E_2^*\rightarrow X$ defines a complete intersection.
Then, for any test form $\varphi$ we have
\[
\int
u^f_{\epsilon_1}\wedge\nabla u^g_{\epsilon_2}\wedge \varphi \rightarrow
(U^f-U^f\wedge R^g).\varphi
\]
as $\epsilon_1,\epsilon_2\rightarrow 0^+$, and the integral to the left belongs
to some Hölder class independently of $\varphi$.
\end{thm}

\begin{proof}
We first note that 
\[
\nabla u^g_{\epsilon_2}=1-\epsilon_2\sum_{\ell\geq 1}
\frac{(\debar s_g)^{\ell-1}}{(|g|^2+\epsilon_2)^\ell},
\]
see the proof of Corollary \ref{sture}. As $U^f\wedge R^f$ is defined as 
the value at zero of the analytic continuation (in the sense of currents) of
$|f|^{2\lambda}u^f\wedge \debar |g|^{2\lambda}\wedge u^g$, what we have to 
prove is 
\begin{eqnarray}\label{eqf}
\lefteqn{\int \frac{s_f\wedge (\debar s_f)^{k-1}}{(|f|^2+\epsilon_1)^k}\wedge
\epsilon_2\frac{(\debar s_g)^{\ell-1}}{(|g|^2+\epsilon_2)^{\ell}}\wedge \varphi
\rightarrow}\\
& &
\int |f|^{2\lambda}u^f_{k,k-1}\wedge \debar |g|^{2\lambda}\wedge u^g_{\ell-1,\ell-2}
\wedge \varphi\Big|_{\lambda=0}\nonumber
\end{eqnarray}
and that the integral on the left belongs to some Hölder class. 
We first consider the case $\ell =1$. The right hand side of \eqref{eqf} should then be 
interpreted as zero. We write the 
integrand on the left hand side of \eqref{eqf} as
$\chi_1(|f|^2/\epsilon_1)\chi_2(|g|^2/\epsilon_2)u^f_{k,k-1}\wedge\varphi$ where 
$\chi_1(t)=t^k/(t+1)^k$ and $\chi_2(t)=1/(t+1)$. 
As in the proof of Theorem \ref{korv} we may assume that 
$u^f_{k,k-1}=\tilde{u}^f_{k,k-1}/\zeta^{k\alpha}$, where $\tilde{u}^f_{k,k-1}$ is a smooth
form, that $|f|^2=|\zeta^{\alpha}|\Phi$ and that $|g|^2=|\zeta^{\beta}|^2\Psi$, where 
$\Phi$ and $\Psi$ are strictly positive 
smooth functions.
Since $\chi_2(\infty)=0$ 
the left hand side of \eqref{eqf} tends to zero and belongs to some Hölder class
by Proposition \ref{tekniskprop}. For $\ell\geq 2$ we proceed as in the proof of Theorem
\ref{satsen} and we see that we may assume that
$f=(f_1,\ldots,f_m)$ and $g=(g_1,\ldots,g_{m_2})$ with 
$f_j=\zeta^{\alpha^j}f'_j$ and $g_j=\zeta^{\beta^j}g'_j$ where all 
$f'_j$ and $g'_j$ are 
non-vanishing and moreover, that for some indices 
$\nu_1$ and $\nu_2$ it holds that
$\zeta^{\alpha}:=\zeta^{\alpha^{\nu_1}}$
divides all $\zeta^{\alpha^j}$ and $\zeta^{\beta}:=\zeta^{\beta^{\nu_2}}$
divides all $\zeta^{\beta^j}$. From the same proof we also see that 
we may assume that $d\bar{\zeta}_j/\bar{\zeta}_j\wedge \varphi$ is smooth
(and compactly supported) for all $\zeta_j$ which divide both $\zeta^{\alpha}$
and $\zeta^{\beta}$, since $f\oplus g$ defines a complete intersection.
We use the notation from the proof of Theorem \ref{satsen}, e.g.\ 
$|f|^2=|\zeta^{\alpha}|^2\Phi=|\zeta^{\alpha'+\alpha''}|^2\Phi$,
$u^f_{k,k-1}=\tilde{u}^f_{k,k-1}/\zeta^{k(\alpha'+\alpha'')}$ and
$|g|^2=|\zeta^{\beta}|^2\Psi=|\zeta^{\beta'+\beta''}|^2\Psi$
etc.
We also introduce the notation $\chi_j(t)$ for the function $t^j/(t+1)^j$, and so,
in particular, we can write $1/(t+\epsilon)^j=\chi_j(t/\epsilon)/t^j$.
For $\ell\geq 2$, one can verify that
\begin{eqnarray}\label{eqg}
\epsilon_2\frac{(\debar s_g)^{\ell-1}}{(|g|^2+\epsilon_2)^\ell}&=&
\frac{1}{\zeta^{(\ell-1)\beta}}\debar\chi_{\ell-1}(|\zeta^{\beta}|^2\Psi/\epsilon_2)
\wedge \tilde{u}^g_{\ell-1,\ell-2}\\
&+&
\frac{1}{\zeta^{(\ell-1)\beta}}
\chi'_{\ell-1}(|\zeta^{\beta}|^2\Psi/\epsilon_2)
\frac{|\zeta^{\beta}|^2}{\epsilon_2}
\frac{\Psi}{\ell-1}
\debar \tilde{u}^g_{\ell-1,\ell-2}. \nonumber
\end{eqnarray}
Using this identity we see that the integral on the left hand side of \eqref{eqf}
splits into two integrals. The integral
corresponding to the last term in \eqref{eqg} tends to zero as 
$\epsilon_1,\epsilon_2 \rightarrow 0$ and belongs to some Hölder class 
according to Corollary \ref{korret}. By Proposition \ref{palle},
the integral corresponding to the first
term on the right hand side of \eqref{eqg} also belongs to some Hölder class 
and tends to
\begin{equation}
-\Big[\frac{1}{\zeta^{k\alpha+(\ell-1)\beta''}}\Big]\otimes
\debar\Big[\frac{1}{\zeta^{(\ell-1)\beta'}}\Big].
\tilde{u}^f_{k,k-1}\wedge\tilde{u}^g_{\ell-1,\ell-2}\wedge\varphi
\end{equation} 
as $\epsilon_1,\epsilon_2 \rightarrow 0$. This is seen to be equal to 
the right hand side of \eqref{eqf} by using the methods in \cite{W}.
\end{proof}

\section{The Passare-Tsikh example}\label{motexsekt}
Let $f=z_1^4$, $g=z_1^2+z_2^2+z_1^3$ and 
$\varphi=\rho\bar{z}_2gdz_1\wedge dz_2$ where $\rho$ has compact support and 
is identically $1$ in a neighborhood of the origin. Since the common zero set
of $f$ and $g$ is just the origin they
define a complete intersection. In \cite{PTmotex} Passare and Tsikh
show that the residue integral  
\[
(\epsilon_1,\epsilon_2)\mapsto I_{f,g}^{\varphi}(\epsilon_1,\epsilon_2)=
\int_{\stackrel{\scriptstyle|f|^2=\epsilon_1}{|g|^2=\epsilon_2}}
\frac{\varphi}{fg}
\]
is discontinuous at the origin. More precisely, they show that for any fixed
positive number $c\neq 1$ one has 
$\lim_{\epsilon\rightarrow 0}I_{f,g}^{\varphi}(\epsilon^4,c\epsilon^2)=0$ but 
$\lim_{\epsilon\rightarrow 0}I_{f,g}^{\varphi}(\epsilon^4,\epsilon^2)\neq 0$.
On the other hand, by Fubini's theorem we have
\begin{eqnarray}\label{eqm}
\int\limits_{[0,\infty)^2}\frac{\epsilon_2\epsilon_2I_{f,g}^{\varphi}(t_1,t_2)dt_1dt_2}
{(t_1+\epsilon_1)^2(t_2+\epsilon_2)^2}&=&
\int \frac{\epsilon_1d|f|^2}{(|f|^2+\epsilon_1)^2}\wedge
\frac{\epsilon_2d|g|^2}{(|g|^2+\epsilon_2)^2}\wedge
\frac{\varphi}{fg}=\nonumber \\
& &
\int \debar \frac{\bar{f}}{|f|^2+\epsilon_1}\wedge
\debar \frac{\bar{g}}{|g|^2+\epsilon_2}\wedge\varphi.
\end{eqnarray}
Hence, this average of the residue integral is continuous at the origin
by Corollary \ref{CR}. In this section we will examine the last integral in
\eqref{eqm} as 
$\epsilon_1,\epsilon_2\rightarrow 0$ explicitly. We will see that it is
continuous at the origin with Hölder exponent at least $1/8$ and that it
tends to zero.
Morally, the value of $I_{f,g}^{\varphi}(\epsilon_1,\epsilon_2)$ at $0$ should be
the Coleff-Herrera current associated to $f$ and $g$ multiplied by
$\bar{z}_2g$ acting on $\rho dz_1\wedge dz_2$. But both $g$ and $\bar{z}_2$ annihilate 
the Coleff-Herrera current since $g$ belongs to 
the ideal generated by $f$ and $g$, and $z_2$ belongs to the 
radical of this ideal. We will thus verify Corollary \ref{CR} 
explicitly in this special case.

Our first objective is to resolve singularities to obtain normal crossings.
This is accomplished by a blow-up of the origin. The map
$\pi\colon\mathcal{B}_0\C^2 \rightarrow \C^2$ looks like 
$\pi(u,v)=(u,uv)$ and $\pi(u',v')=(u'v',u')$ in the two standard coordinate
systems on $\mathcal{B}_0\C^2$. The exceptional divisor, $E$, corresponds
to the sets $\{u=0\}$ and $\{u'=0\}$ and $\pi$ is a biholomorphism
$\mathcal{B}_0\C^2\setminus E \rightarrow \C^2\setminus \{0\}$. 
In the $(u,v)$-coordinates
we have $\pi^*f=u^4$ and $\pi^*g=u^2(1+v^2+u)$. The function $1+v^2+u$
has non-zero differential and its zero locus intersects $E$ normally 
in the two points
$v=i$ and $v=-i$. Moreover, in the $(u',v')$-coordinates we have
$\pi^*f=u'^4v'^4$ and $\pi^*g=u'^2(v'^2+1+u'v'^3)$. The zero locus
of $v'^2+1+u'v'^3$ intersects $E$ normally 
in the points $v'=-i$ and $v'=i$, which we 
already knew, and it does not intersect $v'=0$. 
Also, the differential of $v'^2+1+u'v'^3$ is non-zero on the
zero locus of $v'^2+1+u'v'^3$. 
Hence, $\{\pi^*f\cdot\pi^*g=0\}$ has normal crossings.
We assume that $\varphi$ has support so close to the origin that 
$\mbox{supp}(\pi^*\varphi)\cap\{1+v^2+u=0\}$ has two (compact)
components, $K_1$ and $K_2$, and that these components together with the compacts
$K_3=\mbox{supp}(\pi^*\varphi)\cap\{v=0\}$ and
$K_4=\mbox{supp}(\pi^*\varphi)\cap\{v='0\}$ are pairwise disjoint.
We can then choose a partition of unity $\{\rho_j\}_1^4$ such that
$\sum \rho_j\equiv 1$ on the support of $\pi^*\varphi$ and for each $j=1,2,3,4$,
the support of $\rho_j$ intersects only one of the 
compacts $K_1$, $K_2$, $K_3$ and $K_4$. We choose the numbering such that 
the support of $\rho_j$ intersects $K_j$.
The last integral in \eqref{eqm} now equals
\begin{equation}\label{eqo}
\sum_1^4 \int \debar\frac{\pi^*\bar{f}}{|\pi^*f|^2+\epsilon_1}\wedge
\debar \frac{\pi^*\bar{g}}{|\pi^*g|^2+\epsilon_2}\wedge\rho_j\pi^*\varphi:=
I_1+I_2+I_3+I_4.
\end{equation}
In fact, it is only in $I_3$ we have resonance and we start by considering the 
easier integrals $I_1$, $I_2$ and $I_4$.
The integrals $I_1$ and $I_2$ are similar and we only consider 
$I_1$. The support of $\rho_1$ is contained in a neighborhood of
$p_1=(0,i)$ in the $(u,v)$-coordinates and 
$\rho_1\pi^*\varphi = \rho_1 \pi^*\rho \bar{u}\bar{v}\pi^*g udu\wedge dv$.
Integrating by parts we thus see that
\[
I_1=-\int
\debar\frac{\pi^*\bar{f}}{|\pi^*f|^2+\epsilon_1}
\frac{|\pi^*g|^2}{|\pi^*g|^2+\epsilon_2}\wedge
u\debar ( \bar{u}\bar{v}\rho_1 \pi^*\rho du\wedge dv).
\]
Since $\pi^*f=u^4$ depends on $u$ only, the term of 
$\debar (\bar{u}\bar{v}\rho_1 \pi^*\rho)$ involving $d\bar{u}$ does not give
any contribution to $I_1$. Hence we can replace 
$\debar (\bar{u}\bar{v}\rho_1 \pi^*\rho)$ by $\bar{u}\varphi_1$ where
$\varphi_1$ is smooth and supported where $\rho_1$ is. We put $\zeta_1=u$ and 
$\zeta_2=1+v^2+u$, which defines a change of variables on the support
of $\rho_1$. In these coordinates $\pi^*f=\zeta_1^4$ and 
$\pi^*g=\zeta_1^2\zeta_2$ and so we get
\[
I_1=-\int \frac{1}{\zeta_1^3}\debar \chi(|\zeta_1^4|^2/\epsilon_1) 
\chi(|\zeta_1^2\zeta_2|^2/\epsilon_2)\wedge \bar{\zeta}_1\varphi_1
\]
where $\chi(t)=t/(t+1)$. We also write 
$\debar \chi(|\zeta_1^4|^2/\epsilon_1)=4\tilde{\chi}(|\zeta_1^4|^2/\epsilon_1)
d\bar{\zeta}_1/\bar{\zeta}_1$, where $\tilde{\chi}(t)=t/(t+1)^2$.
To proceed we replace (the coefficient function of)
$d\bar{\zeta}_1/\bar{\zeta}_1\wedge \bar{\zeta}_1\varphi_1$ by its Taylor 
expansion of order one, considered as a function of $\zeta_1$ only, plus a
remainder term $|\zeta_1|^2B(\zeta)$, with $B$ bounded. The terms corresponding
to the Taylor expansion do not give any contribution to $I_1$ since we have 
anti-symmetry with respect to $\zeta_1$ for these terms. Hence, we obtain
\begin{equation}\label{eq001}
|I_1|\lesssim \int_{\Delta}\big|\frac{|\zeta_1|^2B(\zeta)}{\zeta_1^3}
\tilde{\chi}(|\zeta_1^4|^2/\epsilon_1)
\chi(|\zeta_1^2\zeta_2|^2/\epsilon_2)\big|, 
\end{equation}
where $\Delta$ is a polydisc containing the support of $\varphi_1$.
We estimate $|B(\zeta)|$ and $\chi(|\zeta_1^2\zeta_2|^2/\epsilon_2)$ by 
constants, and on the sets 
$\Delta_{\epsilon}=\{\zeta\in \Delta;\, |\zeta_1^4|^2\geq \epsilon_1\}$
and $\Delta\setminus \Delta_{\epsilon}$ we use that 
$\tilde{\chi}(|\zeta_1^4|^2/\epsilon_1)\lesssim \epsilon_1/|\zeta_1^4|^2$
and 
$\tilde{\chi}(|\zeta_1^4|^2/\epsilon_1)\lesssim |\zeta_1^4|^2/\epsilon_1$
respectively, to see that the right hand side of \eqref{eq001} is of the size
$|\epsilon|^{1/8}$.

To deal with $I_4$ we proceed as follows. The support of $\rho_4$ is 
contained in a neighborhood of $p_4=(0,0)$ in the $(u',v')$-coordinates and 
$\pi^*f=u'^4v'^4$ and $\pi^*g=u'^2(1+v'^2+u'v'^3):=u'^2\tilde{g}$. On the support of 
$\rho_4$ we have $\tilde{g}\neq 0$. 
The multiindices $(4,4)$ and $(2,0)$ are linearly independent
and so we can make the factor $\tilde{g}$ disappear. Explicitly, choose a square root 
$\tilde{g}^{1/2}$ of $\tilde{g}$ and put $\zeta_1=u'\tilde{g}^{1/2}$ and
$\zeta_2=v'\tilde{g}^{-1/2}$. In these coordinates $\pi^*f=\zeta_1^4\zeta_2^4$
and $\pi^*g=\zeta_1^2$. One also checks that 
$\rho_4\pi^*\varphi=|\zeta_1|^2\pi^*g \varphi_4$ where $\varphi_4$ is a test form of
bidegree $(2,0)$. After an integration by parts we see that
\begin{equation}\label{eq002}
I_4=\int \frac{\pi^*\bar{f}}{|\pi^*f|^2+\epsilon_1}
\debar \frac{|\pi^*g|^2}{|\pi^*g|^2+\epsilon_2}\wedge
\debar(|\zeta_1|^2\varphi_4).
\end{equation}
Since $\pi^*g=\zeta_1^2$ only depends on $\zeta_1$ we may replace 
$\debar(|\zeta_1|^2\varphi_4)$ by $|\zeta_1|^2\debar \varphi_4$ in \eqref{eq002}.
Computing $\debar (|\pi^*g|^2/(|\pi^*g|^2+\epsilon_2))$ we find that
\[
I_4=2\int \frac{1}{\zeta_1^3\zeta_2^4}
\chi(|\zeta_1^4\zeta_2^4|^2/\epsilon_1)
\tilde{\chi}(|\zeta_1^2|^2/\epsilon_2)d\bar{\zeta}_1\wedge \debar \varphi_4.
\]
With abuse of notation we write the test form $d\bar{\zeta}_1\wedge \debar \varphi_4$
as $\varphi_4d\zeta\wedge d\bar{\zeta}$.
Let $M=M_{1,2}^{1,2}$ be the operator defined in Lemma \ref{taylorlemma}.
Explicitly, we have 
\begin{eqnarray*}
M\varphi_4&=&
M_1^1\varphi_4+M_2^2\varphi_4-M_1^1M_2^2\varphi_4\\
&=&
M_1^1(\varphi_4-M_2^2\varphi_4)+M_2^2(\varphi_4-M_1^1\varphi_4)+
M_1^1M_2^2\varphi_4.
\end{eqnarray*}
All of the following properties will not be important for this computation but to 
illustrate Lemma \ref{taylorlemma} we note that
the second expression of $M\varphi$ reveals that $M\varphi_4$ can be written as a sum 
of terms $\phi_{IJ}(\zeta)\zeta^I\bar{\zeta}^J$ with $I_1+J_1\leq 1$ and 
$I_2+J_2\leq 2$ and moreover, that $\phi_{IJ}$ is independent of at least
one variable and is of the size $\mathcal{O}(|\zeta_1|^2)$ if it depends on 
$\zeta_1$ and of the size
$\mathcal{O}(|\zeta_2|^3)$ if it depends on $\zeta_2$.
By Lemma \ref{taylorlemma} we also have 
$\varphi_4=M\varphi_4+|\zeta_1|^2|\zeta_2|^3B(\zeta)$ for some bounded function
$B$ and so 
\[
I_4=\int_{\Delta}
\frac{1}{\zeta_1^3\zeta_2^4}\chi\tilde{\chi}M\varphi_4+\int_{\Delta}
\frac{1}{\zeta_1^3\zeta_2^4}\chi\tilde{\chi}|\zeta_1|^2|\zeta_2|^3B(\zeta)=:
I_{4.1}+I_{4.2},
\]
where $\Delta$ is a polydisc containing the support of $\varphi_4$. By anti-symmetry
$I_{4.1}=0$. 
To estimate $I_{4.2}$ we use that $|\chi B|$ is bounded by a constant and that 
$\tilde{\chi}(\Psi |\zeta_1^2|^2/\epsilon_2)\lesssim \epsilon_2/|\zeta_1^2|^2$ 
and 
$\tilde{\chi}(\Psi |\zeta_1^2|^2/\epsilon_2)\lesssim |\zeta_1^2|^2/\epsilon_2$
on the sets 
$\Delta_{\epsilon}=\{\zeta\in \Delta;\, |\zeta_1^2|^2\geq \epsilon_2\}$ and 
$\Delta\setminus \Delta_{\epsilon}$ respectively. Hence,
\begin{equation}\label{eq006}
|I_{4.2}|\lesssim \int_{\Delta_{\epsilon}}
\frac{\epsilon_2}{|\zeta_1^2|^2|\zeta_1||\zeta_2|}+
\int_{\Delta\setminus\Delta_{\epsilon}}
\frac{|\zeta_1^2|^2}{\epsilon_2|\zeta_1||\zeta_2|},
\end{equation}
which is seen to be of the size $|\epsilon|^{1/4}$.

It remains to take care of $I_3$. We are now working close to $u=v=0$ and 
$\pi^*f=u^4$ and $g=u^2(1+v^2+u):=u^2\tilde{g}$. The multiindices are linearly dependent
and we cannot dispose of the non-zero factor $\tilde{g}$. We rename our variables 
$(u,v)=(\zeta_1,\zeta_2)$ and proceed in precisely the same way as we did when we were
considering $I_1$. We get
\[
I_3=-4\int \frac{1}{\zeta_1^3}\tilde{\chi}(|\zeta_1^4|^2/\epsilon_1) 
\chi(\Phi|\zeta_1^2|^2/\epsilon_2) \varphi_3
d\zeta\wedge d\bar{\zeta},
\]
where $\Phi=|\tilde{g}|^2$ is a strictly positive smooth function and $\varphi_3$ is 
smooth with compact support. As before, we replace $\varphi_3$ by 
$M_{\zeta_1}^1\varphi_3+|\zeta_1|^2B(\zeta)$. The integral corresponding
to $|\zeta_1|^2B(\zeta)$ satisfies the same estimate as the one in \eqref{eq001} and hence
is of the size $|\epsilon_1|^{1/8}$. We cannot
use anti-symmetry directly to conclude the the integrals corresponding to the other
terms in the Taylor expansion tend to zero since the factor $\tilde{g}$ is present. 
We illustrate why this is true anyway by considering the integral corresponding to the
term $\varphi_3(0,\zeta_2)$. Let $\Delta$ be a polydisc containing the support of 
$\varphi_3$ and consider 
\begin{equation}\label{tom}
\int_{\Delta} \frac{1}{\zeta_1^3}\tilde{\chi}(|\zeta_1^4|^2/\epsilon_1) 
\chi(\Phi|\zeta_1^2|^2/\epsilon_2)\varphi_3(0,\zeta_2).
\end{equation}
We introduce the smoothing parameter 
$t=|\zeta_1^2|^2/\epsilon_2$ as an independent variable and write
\[
\chi(\Phi t)=\chi(\Phi t)-M_{\zeta_1}^1\chi(\Phi t)+M_{\zeta_1}^1\chi(\Phi t):=
|\zeta_1|^2B(t,\zeta)+M_{\zeta_1}^1\chi(\Phi t).
\]
Here $B$ is bounded on $[0,\infty]\times \Delta$. Substituting into \eqref{tom} we obtain
one integral corresponding to $|\zeta_1|^2B(|\zeta_1^2|^2/\epsilon_2,\zeta)$, 
which satisfies an estimate like
\eqref{eq001}, while the integral corresponding to 
$M_{\zeta_1}^1\chi(\Phi |\zeta_1^2|^2/\epsilon_2)$ is zero
since we have anti-symmetry with respect to $\zeta_1$.
Hence $|I_3|\lesssim |\epsilon|^{1/8}$.

\bigskip

{\bf Acknowledgments}

\medskip

I am grateful to Mats Andersson and Jan-Erik Bj\"ork for valuable comments and remarks 
on preliminary versions.

\bibliography{residyartikel}

\end{document}